\documentclass[11pt,leqno]{article}

\usepackage{graphicx}
\usepackage{latexsym} 
\usepackage{amsmath}
\usepackage{amssymb}
\usepackage{amsfonts}
\usepackage{enumerate}
\usepackage{theorem}

\usepackage{tikz}

\usepackage{xcolor}

\theoremstyle{plain}
\newtheorem{teo}{Theorem}[section]
\newtheorem{lem}[teo]{Lemma}
\newtheorem{cor}[teo]{Corollary}

{\theorembodyfont{\rmfamily}
\newtheorem{defin}[teo]{Definition}
\newtheorem{oss}[teo]{Remark}

}

\renewcommand{\eqref}[1]{\textnormal{(\ref{#1})}}

\numberwithin{equation}{section}

\newcommand{\cvd}{\hfill$\square$}
\newcommand{\proof}[1]{\noindent\textsc{Proof#1}}
\newcommand{\rmi}{\mathrm{i}}
\newcommand{\rme}{\mathrm{e}}

\title{Stable determination of a rigid scatterer in
elastodynamics}

\author{Luca Rondi\thanks{Dipartimento di Matematica,
Universit\`a di Milano, Italy. E-mail: \texttt{luca.rondi@unimi.it}} \and\
Eva Sincich\thanks{Dipartimento di Matematica e Geoscienze, Universit\`a di Trieste, Italy. E-mail: \texttt{esincich@units.it}} \and\
Mourad Sini\thanks{RICAM, Austrian Academy of Sciences, Austria. E-mail: \texttt{mourad.sini@oeaw.ac.at}}
}

\date{}

\begin{document}

\maketitle

\setcounter{section}{0}
\setcounter{secnumdepth}{2}

\begin{abstract}
\noindent
We deal with an inverse elastic scattering problem for the shape determination of a rigid scatterer in the time-harmonic regime. We prove a local stability estimate of $log\,log$ type for the identification of a scatterer by a single far-field measurement. The needed a priori condition on the closeness of the scatterers is estimated by the universal constant appearing in the Friedrichs inequality.  

\smallskip

\noindent\textbf{AMS 2020 Mathematics Subject Classification}
Primary 35R30. Secondary 74B05, 35P25.


\smallskip

\noindent \textbf{Keywords} inverse problems, scattering, linear elasticity, elastic waves, stability
\end{abstract}

\section{Introduction}\label{intro}

We consider the scattering of time-harmonic elastic waves by a rigid scatterer in $\mathbb{R}^N$ with $N\ge 2$. 
The time-harmonic elastic waves in a homogeneous and isotropic elastic medium satisfy the Navier equation
\begin{equation}\label{Naviereq-INTRO}
\mu\Delta u+(\lambda+\mu){\nabla}(\mathrm{div} (u))+\rho \omega^2 u=0,
\end{equation}
where $\lambda$ and $\mu$ are the Lam\'e constants such that $\mu>0$ and $\lambda+2\mu>0$, $\rho>0$ is the density and $\omega>0$ is the frequency. By the Helmholtz decomposition, any solution $u$ to \eqref{Naviereq-INTRO} is the superposition of a longitudinal wave $u_p$ and a transversal wave $u_s$, which are solutions to the Helmholtz equation with wave numbers $\omega_p=\sqrt{\frac{\rho}{\lambda + 2\mu}}\;\omega$ and $\omega_s=\sqrt{\frac{\rho}{\mu}}\;\omega$, respectively.

If an incident wave $u^{inc}$, which is usually given by  an entire solution to \eqref{Naviereq-INTRO}, meets a rigid scatterer $K$, then it is perturbed by the formation of a scattered wave $u^{scat}$ outside $K$. The total field $u$ is the superposition of the incident and the scattered wave and, for a rigid scatterer, satisfy the following Dirichlet boundary condition on the boundary of the scatterer
\begin{equation}\label{BC-INTRO}
u=0 \qquad \text{on }\ \partial K.
\end{equation}
The scattered wave $u^{scat}$ is characterised by being a radiating solution to \eqref{Naviereq-INTRO}, namely
its longitudinal wave $u^{scat}_p$ and transversal wave $ u^{scat}_s$ are radiating solutions to the corresponding Helmholtz equations. The radiation condition for elastic waves is usually referred to as the Kupradze radiation condition.  

As incident wave $u^{inc}$, we take either a \emph{longitudinal plane wave} 
\begin{equation}\label{planarlong-INTRO}
u^{inc}_p(x)=d\;\rme^{\rmi \omega_pd\cdot x},\qquad x\in \mathbb{R}^N,
\end{equation}
where $d\in \mathbb{S}^{N-1}$ is the \emph{direction of incidence}, or a \emph{transversal plane wave}
\begin{equation}\label{planartrans-INTRO}
u^{inc}_s(x)=p\; \rme^{\rmi \omega_sd\cdot x},\qquad x\in \mathbb{R}^N,
\end{equation}
where $p\in \mathbb{C}^N\backslash \{0\}$ is a unitary vector orthogonal to $d$.
We can also consider a linear combination of longitudinal and transversal plane waves, namely
\begin{equation}\label{planarcomb-INTRO}
u^{inc}(x)=c_pu^{inc}_p(x)+c_su^{inc}_s(x),\qquad x\in \mathbb{R}^N
\end{equation}
for some $c_p,c_s\in\mathbb{C}$ such that $|c_p|^2+|c_s|^2=1$, in such a way to have
\begin{equation}\label{planarcomb-prop-INTRO}
\|u^{inc}(x)\|=1,\qquad x\in \mathbb{R}^N
\end{equation}

The forward scattering problem for a rigid obstacle is classical and, under mild regularity assumptions on the obstacle $K$, it is well-known to have a unique solution.

By the  Kupradze radiation condition, the scattered wave $u^{scat}$ has the following asymptotic behaviour 
 \begin{equation}\label{decay+far-field INTRO}
u^{scat}(x; d)=
\frac{\rme^{\rmi \omega_pr}}{r^{(N-1)/2}}
U_p(\hat{x}; d)+\frac{\rme^{\rmi \omega_s r}}{r^{(N-1)/2}}
U_s(\hat{x}; d)
+O\left(\frac{1}{r^{(N+1)/2}}\right),
\end{equation}
\noindent
as $r=\|x\|$ goes to $+\infty$, uniformly in all directions $\hat{x}=x/\|x\|\in \mathbb{S}^{N-1}$.
The vector fields $U_p$ and $U_s$ are called longitudinal and transversal far-field patterns, respectively. Since they characterise, respectively, the asymptotic behaviour of 
the normal and of the tangential component, with respect to $\mathbb{S}^{N-1}$, of $u^{scat}$, by measuring the asymptotic behaviour of $u$, or equivalently of $u^{scat}$, as $r$ goes to $+\infty$, both the longitudinal part and the transversal part of the far-field pattern of $u^{scat}$ can be measured. 

We are concerned with the following geometrical inverse problem in the context of linear elasticity. Given an incident wave $u^{inc}$, one can measure the vector fields $(U_p(\cdot, d), U_s(\cdot, d))$, which as usually referred to as the corresponding scattering data. By changing the incident wave, for instance by changing the frequency $\omega$ or the incident direction $d$, one can obtain different scattering data.
We wish to determine the scatterer $K$ by using as measured data the scattering data corresponding to one or more incident waves.

The unique determination of $K$ using the measured data corresponding to all the incident directions $d\in \mathbb{S}^{N-1}$, with a fixed frequency $\omega$, was first shown in \cite{Haen-Hsia}. In their work, they use both the components $U_p(\cdot, d)$ and $U_s(\cdot, d)$ of the elastic farfields. Next, it was proved that actually only one component of the farfield $U_p(\cdot, d)$ or $U_s(\cdot, d)$ is enough, meaning that either the pressure or the shear waves are enough to uniquely determine the scatterer $K$. This result was justified first for $C^4$-smooth scatterers in \cite{Gi-Si} and later it was extended to Lipschitz-smooth ones in \cite{Kar-Si}. In addition, reconstruction schemes were proposed in \cite{Hu-Ki-Si, Kar-Si-2} to actually reconstruct the scatterer $K$.    

Here, we are interested in the determination of the scatterer $K$ by the knowledge of the longitudinal and transversal far-field patterns corresponding to a single incident wave provided some suitable a priori information  on the location of the scatterer is known.  

A special instance of such a problem has been previously analysed in \cite{Gi-Mi} in a two dimensional setting. Indeed the authors proved a uniqueness result if it is a priori known that possible  scatterers do not deviate too much in area, or more precisely under the following closeness condition 
 \begin{eqnarray}\label{closeGM}
 |K\Delta K'|\le \frac{ k^2_{0,1}\pi \mu}{\rho \omega^2}
 \end{eqnarray}
where
$k_{0,1}\approx 2.4048$ is the first zero of the Bessel function $J_0$. Their arguments is strongly based on the fact that a lower estimate for the first Dirichlet eigenvalue of the negative Lam\'e operator in $K\Delta K'$ in terms of its Lebesgue measure can be achieved by the use of the  Faber-Krahn inequality (see also \cite{Kaw}).

Here we study the stability issue for the same problem in any dimension $N\ge 2$. We prove a $log\,log$ type stability estimate for the unknown scatterer under a slightly stronger a priori closeness condition (see Section~\ref{mainresultsec} for a precise statement), namely we assume that $K$ and $K'$ are both contained in a given scattered $K^+$ and that
 \begin{eqnarray}\label{closeRSS}
 |K^+\backslash (K\cap K')|< H_1=\left(\frac{\min\{2\mu,2\mu+\lambda\}}{64 C(N)^2\rho\omega^2}\right)^{N/2}
 \end{eqnarray}
 where $C(N)$ is an absolute constant depending on the dimension $N$ only, actually it is the one of the isoperimetric inequality, see \eqref{constant}. Just for comparison, for $N=2$ and assuming for simplicity $\lambda>0$, our closeness bound becomes
$$ H_1= \frac{\pi\mu}{8\rho\omega^2}.$$
Such a slightly more restrictive a priori bound is justified by the fact that, in order to deal with stability, we are led to replace the use of the Faber-Krahn inequality with the one of the Friedrichs inequality. Besides the closeness condition, we require some a priori regularity of the unknown scatterer, in particular we require it to be of class $C^{2,\alpha}$, $0<\alpha<1$. We note, however, that we allow $K$ to have more than one connected component. Our stability estimate, which is the main result of the paper, is stated in Theorem~\ref{mainthm}.

Let us note that even if the estimate is rather weak, being of $log\, log$ type, this is rather common for these kinds of inverse scattering or boundary value problems. Moreover, it has been shown that a single $log$ estimate is optimal for the stability of these inverse problems even if many measurements are performed, see \cite{Man} and \cite{DC-R}.

We recall that analogous local uniqueness  and local stability results have been previously achieved in the acoustic framework in \cite{Gi,St-U} and \cite{S-S}, respectively, by means of a spectral type approach as then extended to elasticity in \cite{Gi-Mi}. 
Unfortunately, as already observed, these arguments cannot be applied to extend the stability result in the elastic case and hence new tools and an original strategy have to be introduced. 

Stability results of $log\,log$ type related to the elasticity system are derived in \cite{Hig, Mor-Ros-ip, Mor-Ros-mem} in the stationary case, that is, $\omega=0$, with a single pair of displacement and traction fields measured on a surface surrounding the unknown scatterer $K$. Both rigid inclusions and cavities have been treated. More recently, \cite{Mor-Ros-Ves2}, an optimal single $log$ estimate has been obtained for the determination of cavities in the two dimensional case, by exploiting an optimal three-spheres inequality at the boundary, \cite{Al-Ros-Ves}, or a doubling inequality at the boundary, \cite{Mor-Ros-Ves}, that have been established for the Kirchhoff-Love plate's equation. These kinds of inequality might allow us to improve our estimate to a single $log$ one as well, but to the best of our knowledge none of them is already available in our setting.

Among other results of identification of targets by means of a single or finitely many far-field data under a priori geometric constraints, we wish to recall the well-known uniqueness result in \cite{CoSl} for small obstacles in an acoustic context and the corresponding stability estimates in \cite{Isak92,Isak93}, under the additional starshapedness hypothesis of the scatterer.

Another case in which one measurement, or at least few measurements, uniquely identifies a scatterer is when the scatterer
satisfies a different strong geometric condition, namely it is assumed to be polyhedral. In the case of obstacles, this means that the objects to be determined are (a collection of) polygons in dimension 2 or polyhedra in dimension 3. In this direction, in \cite{C-Y} it was proved the first uniqueness result for sound-soft scatterers in the acoustic framework. In the same case, an optimal uniqueness result with a single measurement was established in \cite{AR}. These results have been extended to a variety of other boundary conditions and to the electromagnetic case by several authors. The first stability result,
still for sound-soft scatterers in the acoustic case, was obtained in \cite{Ron08}, and it was followed by analogous results for the sound-hard boundary condition \cite{LPRX} and for the electromagnetic case \cite{L-R-X}.
This line of research has been first extended to the elastic case in \cite{El1}, where the third and fourth boundary conditions are considered and uniqueness is established with two (suitable) measurements or one (suitable) measurement, respectively. For polyhedral obstacles, one (suitable) measurement is enough to determine the obstacle and the boundary condition, provided the latter is still of the third or fourth type, see \cite{L-X} where a corresponding stability estimate is also proved. In \cite{El2} it is shown that one measurement is enough to uniquely determine a rigid polygon. More recently, still in dimension 2, in \cite{Diao-et-al} it is showed that four measurements allow to uniquely recover a collection of polygons as well as their mixed boundary condition, since the obstacles may, at the same time, be purely rigid or be traction-free or satisfy an impedance boundary condition on different parts of the boundary.

Finally, we mention that in \cite{H-N-S} uniqueness results using few measurements are derived even if the scatterers
neither satisfy smallness conditions nor have polygonal or polyhedral shapes. Instead, the authors assume that the boundary of the scatterers are nowhere analytic and show that for the Dirichlet boundary condition one measurement is enough while for the Neumann boundary condition $N-1$ measurements are enough (recalling that $N$ is the space dimension). These results are proved for the Helmholtz model for any $N$, and then generalised to a larger family of elliptic second order operators when $N\leq 3$.
      
In order to analyse the inverse problem we preliminary discuss the direct one. Indeed, in Section~\ref{prelsec}, we observe that the direct scattering problem is well-posed (see Theorem~\ref{existenceteo}) recalling a classical result due to Kupradtze et al \cite{Kup2}. In Theorem~\ref{regularityteo} we prove a regularity result up to the boundary $\partial K$ for the solution $u$, independent on the scatterer $K$. The proof, which may be found in the Appendix, is based on well-known regularity estimates for elliptic systems with Dirichlet boundary condition \cite{ADN} and a preliminary bound for the solution obtained by a continuity argument inspired by Mosco convergence, as done in \cite{Men-Ron} for the acoustic case and in \cite{L-R-X} for the electromagnetic one. In particular, we also obtain a uniform decay property, as $r\to+\infty$, of $u^{scat}$, again independent on the scatterer $K$. We conclude this preliminary part by reviewing, in Subsection~\ref{Friedrichs inequality}, the Friedrichs inequality. We observe that its constant depends on the dimension only and can be explicitly evaluated. This allows us to state the closeness condition with an explicit constant depending on the dimension $N$ and the coefficients of the Navier equation only.

In Section~\ref{mainresultsec} we state our main stability result, Theorem~\ref{mainthm}, whose proof is developed in 
Section~\ref{proofsec}.

The strategy of the proof is the following. Assuming we have two scatterers $K$ and $K'$ satisfying the closeness condition,
we wish to estimate their Hausdorff distance $\tilde{d}$ from the difference of the corresponding far-field patterns.
From the error on the far-field patterns, we estimate the error on the total field in a region surrounding the scatterers by a classical far-field to near-field estimate for the Helmholtz equation applied to the longitudinal and transversal part of the scattered wave. This estimate has been proved first in dimension $3$ in \cite{Isak92}, see also \cite{Bus}, and then generalised to any dimension in \cite{R-S}. For a suitably chosen small parameter $s>0$, we call $V_s$ the region outside $K\cup K'$ whose points can be reached from infinity by a suitable tube of radius $s$.
By a standard unique continuation argument, whose main ingredient is a three-spheres inequality for the Helmholtz equation proved in \cite{Bru} which is iteratively applied inside the $s$-tube to 
the longitudinal and transversal part of the field, we are able to estimate the error on the total field on the boundary of $V_s$, see Lemma~\ref{firstcontinuation}.
Up to swapping $K$ with $K'$, we may find $A_s$, a suitable connected component of $\mathbb{R}^N\backslash (K\cup \overline{V_s})$, which contains a ball $B$ of radius proportional to $\tilde{d}$. By the regularity of $K$ and $K'$ and the closeness condition, we infer $A_s$ is a set of finite perimeter (with a perimeter depending on the a regularity of the scatterers only) whose measure is strictly less than the closeness constant $H_1$. For details on the definition of $V_s$ and $A_s$ see page~\pageref{V_sdefin} and following pages.

By the estimate on the error of the field on the boundary of $V_s$ and the boundary condition, we obtain an estimate on the solution $u$ on the boundary of $A_s$. It is now that the Friedrichs inequality comes into play. By using the Friedrichs inequality, the bound on the boundary of $A_s$, the a priori bound on the solution $u$ and the idea of the proof of the first Korn inequality, we are able to estimate the $L^2$-norm of $\nabla u$ in $A_s$ by a quantity depending on $s$ and the norm of $u$ on $\partial A_s$, see Lemma~\ref{cruciallem}. Actually, we can estimate the $L^{\infty}$-norm of $u_p$ and $u_s$ on the ball $B$ by a constant depending on the far-field error only, provided we choose a suitable value of $s$, see 
Corollary~\ref{cor} and Remark~\ref{firststeposs}.

From $B$ we move towards infinity again and, by a unique continuation argument pretty similar to the one used before, 
we show that $u$ has to be small even far away from the scatterers, such a smallness depending on the smallness of
$u_p$ and $u_s$ on the ball $B$ and on the radius of $B$ itself. However, by our choice of the incident field and the decay of the scattered one, the total field $u$ can not be too small far from $K$. Combining these two pieces of information, we are  finally able to estimate the radius of $B$ and consequently the Hausdorff distance between $K$ and $K'$, see Lemma~\ref{gettingoutlemma} and Remark~\ref{anotheross}. This second part of the proof is inspired by an analogous procedure developed in \cite{Ron08}.

By this technique we obtain a stability estimate which is however extremely weak being of $log\, log\, log$ type. This is due to the fact that the estimate of the error on the boundary of $V_s$ is already of $log\, log$ type since $V_s$ can be extremely irregular. Another $log$ comes from the second part of the procedure when we move from $B$ towards infinity.
On the other hand, we can apply a refining procedure which is by now standard, see for instance \cite{Al-etal}. Provided the error on the farfields is small enough, the two scatterer are close enough and, by their a priori regularity, it can be inferred that the unbounded connected component of $\mathbb{R}^N\backslash (K\cup K')$, which we call $V_0$, satisfies a Lipschitz type regularity, see Lemma~\ref{RLGlemma} which is an easy consequence of \cite[Lemma~8.1]{Al-etal}.
By such regularity, we can improve our unique continuation estimate up the boundary of $V_0$ to a single $log$ estimate. By using a suitable domain $A_0$ and by exactly the same procedure as before, we are then able to improve our estimate 
to the final one of $log\, log$ type.

To conclude, we wish to put in evidence a delicate point of the proof that also explains the presence of the scatterer $K^+$ in our closeness condition. Even if $K$ and $K'$ are smooth, the domain $V_0$ can be extremely irregular, unless we know that the scatterer are close enough to apply Lemma~\ref{RLGlemma}. This implies that we are not able to estimate the error on the field up to the boundary of $V_0$ from the farfields error and are thus forced to introduce the set $V_s$, $s>0$. In turn, this introduces another difficulty. In fact the domain $A_0$, which we can construct from $V_0$, is contained in $K'\backslash K$, thus its measure is bounded by the measure of $K\Delta K'$. Instead, for $s>0$, the measure of the domain $A_s$, which we can construct from $V_s$, is not controlled by the measure of $K\Delta K'$. This is the reason why we need to introduce $K^+$ in the closeness condition.

\subsubsection{Acknowledgements} LR is partly supported by GNAMPA--INdAM through Projects 2018 and 2019. The work of  ES is partly performed under the PRIN grant No. 201758MTR2-007. MS is partly supported by the Austrian Science Fund (FWF): P 30756-NBL. Part of this work was done during visits to the University of Trieste, Italy, and to the Radon Institute, Austria. The authors wish to thank both the institutions for their kind hospitality.

\section{Preliminaries}\label{prelsec}
Throughout the paper the integer $N\geq 2$ will denote the space dimension. We note that we drop the dependence of any constant from the space dimension $N$.
 For any two column vectors $U=(U^1,\ldots,U^N)^T$ and $V=(V^1,\ldots,V^N)^T$ in $\mathbb{C}^N$
$$U\cdot V=U^TV=
\sum_{i=1}^N U^iV^i.$$
Here, and in the sequel,
for any matrix $A$, $A^T$ denotes its transpose. By $I_N$ we denote the identity $N\times N$ matrix.
For any two matrices $A=\{a_{i,j}\}_{i,j=1}^N$ and $B=\{b_{i,j}\}_{i,j=1}^N$,
$$A : B= \sum_{i,j=1}^N a_{ij}b_{ij}.$$

For any $x=(x_1,\ldots,x_N)\in\mathbb{R}^N$,  we denote $x=(x',x_N)\in\mathbb{R}^{N-1}\times \mathbb{R}$.
For any $s>0$ and any $x\in\mathbb{R}^N$, 
$B_s(x)$ denotes the open ball contained in $\mathbb{R}^N$ with radius $s$ and centre $x$, whereas $B_s=B_s(0)$.
For any $E\subset \mathbb{R}^N$, $B_s(E)=\bigcup_{x\in E}B_s(x)$.
Given a point $x\in\mathbb{R}^N$, a vector $v\in\mathbb{S}^{N-1}$, and constants
$r>0$ and $\theta$, $0<\theta\leq\pi/2$, we call $\mathcal{C}(x,v,r,\theta)$ the open cone with vertex in $x$, bisecting vector given by $v$, radius $r$ and amplitude given by $\theta$, that is
$$\mathcal{C}(x,v,r,\theta)=\left\{y\in\mathbb{R}^N:\ 0<\|y-x\|<r\text{ and }\cos(\theta)<\frac{y-x}{\|y-x\|}\cdot v\leq 1\right\}.$$
We remark that by a cone we always mean a bounded not empty open cone of the kind defined above.

For any measurable subset of $\mathbb{R}^N$ we call $|E|$ its $N$-dimensional Lebesgue measure. By $\mathcal{H}^{N-1}$ we denote the $(N-1)$-dimensional Hausdorff measure.

\begin{defin}
Let $\Omega\subset\mathbb{R}^N$ be a bounded open set. Let $k$ be a nonnegative integer and $0\leq\alpha\leq 1$.

We say that
$\Omega$ is of 
class $C^{k,\alpha}$ (Lipschitz if $k=0$ and $\alpha=1$, $C^k$ if $\alpha=0$)
if for any $x\in\partial\Omega$ there exist a $C^{k,\alpha}$
function $\phi_x:\mathbb{R}^{N-1}\to\mathbb{R}$ and a neighbourhood $U_x$ of $x$
such that for any $y\in U_x$ we have, up to a rigid transformation depending on $x$,
$$y=(y',y_N)\in\Omega\quad \text{if and only if}\quad  y_N<\phi_x(y').$$

We also say that $\Omega$ is of 
class $C^{k,\alpha}$ (Lipschitz if $k=0$ and $\alpha=1$, $C^k$ if $\alpha=0$) with positive constants $r$ and $L$ if for any $x\in\partial\Omega$ we can choose $U_x=B_r(x)$ and $\phi_x$ such that $\|\phi_x\|_{C^{k,\alpha}(\mathbb{R}^{N-1})}\leq L$.
\end{defin}

\begin{oss}\label{normaloss}
If $k+\alpha>1$, and $\Omega$ is an open set of
class $C^{k,\alpha}$ with constants $r$ and $L$, there exists positive constants $r_1$ and $L_1$, depending on $k$,
$\alpha$, $r$ and $L$ only, such that $\Omega$ is of class $C^{k,\alpha}$ with constants $r_1$ and $L_1$ with the further condition that for any $x\in\partial \Omega$ we have $\nabla \phi_x(x')=0$. Therefore, without loss of generality, whenever $k+\alpha>1$ we tacitly assume that this condition is satisfied all over $\partial\Omega$.
\end{oss}

\subsection{The direct scattering problem}

We say that $\Omega\subset\mathbb{R}^N$ is a \emph{domain} if it is open and connected. We say that $\Omega$ is an \emph{exterior domain} if it is a domain containing the exterior of a ball. We say that $K\subset\mathbb{R}^N$ is a \emph{scatterer} if $K$ is compact and $\Omega=\mathbb{R}^N\backslash K$ is connected, that is, $\Omega$ is an exterior domain. We say that a scatterer $K$ is an \emph{obstacle} if $K=\overline{D}$ where $D$ is an open set which we can pick as the interior of $K$.

We consider the inverse scattering problem for the Navier equation modelling time-harmonic elastic waves in a homogeneous and isotropic elastic medium under the presence of a rigid scatterer. For the direct scattering problem, which we here describe, we refer to the classical works of Kupradze and others
\cite{Kup1,Kup2} and to the more recent one \cite{Bao-et-al}, for instance.
Let us consider, in an open set $\Omega\subset\mathbb{R}^N$, $N\geq 2$, a weak solution $u$ to the \emph{Navier equation}
\begin{equation}\label{Navier}
\mu\Delta u+(\lambda+\mu)\nabla^T(\mathrm{div} (u))+\rho \omega^2 u=0\qquad \text{in }\Omega. 
\end{equation}
Here $\lambda$ and $\mu$ are the Lam\'e constants such that $\mu>0$ and $\lambda+2\mu>0$, $\rho>0$ is the density and 
$\omega>0$ is the frequency. We assume all these parameters to be constants. 
The function $u=(u^1,\ldots,u^N)^T$, the so-called \emph{field} of the time-harmonic wave, is assumed to be a column vector. We note that if $v$ is a scalar function, we often use
$\nabla^T$ to denote the column vector $(\nabla v)^T$.

A vector field
$u\in H^1_{loc}(\Omega,\mathbb{C}^N)$ is a weak solution to \eqref{Navier} if for any $v\in H^1(\Omega,\mathbb{C}^N)$ with compact support in $\Omega$ we have
\begin{equation}\label{Navierweak}
2\mu\int_{\Omega}Eu : \overline{Ev}+\lambda\int_{\Omega}\mathrm{div} (u) \overline{\mathrm{div} (v)} - \rho\omega^2\int_{\Omega}u\cdot \overline{v}=0.
\end{equation}
Here $Eu=\frac{1}{2}(\nabla u+(\nabla u)^T)$ denotes the symmetric gradient of $u$.
Hence, $\nabla u=\{u^i_j\}_{i,j=1}^N$ and $Eu=\frac{1}{2}\{u^i_j+u^j_i\}_{i,j=1}^N$.

In linearised elasticity $Eu$ corresponds to the \emph{strain tensor} and, by Hooke's law, the \emph{stress} $\sigma$ is given by 
$$\sigma(u)=2\mu Eu+\lambda \mathrm{tr}(Eu)I_N=
2\mu Eu+\lambda\mathrm{div}(u)I_N,$$
where $\mathrm{tr}$ denotes the trace.
In particular, we have
\begin{equation}\label{elliptic}
\min\{2\mu,2\mu+\lambda\}(Eu:\overline{Eu})\leq \sigma(u):\overline{\sigma(u)}\leq \max\{2\mu,2\mu+\lambda\}(Eu:\overline{Eu}).
\end{equation}
We note that
\eqref{Navier} can be rewritten as
$$\mathrm{div}(\sigma(u))+\rho \omega^2 u=0\qquad \text{in }\Omega,$$
where the $\mathrm{div}$ applies row by row. In fact,
$$\mathrm{div}(\mathrm{div}(u)I_N)=\mathrm{div}((\nabla u)^T)=\nabla^T(\mathrm{div}(u)).$$

We call
$$\mathcal{K}(\Omega)=\{u\in L^2(\Omega,\mathbb{C}^N):\ Eu\in L^2(\Omega,\mathbb{C}^{N\times N})\},$$
which is a Hilbert space with the corresponding norm
$$\|u\|_{\mathcal{K}(\Omega)}=\left(\|u\|^2_{L^2(\Omega)}+\|Eu\|^2_{L^2(\Omega)}\right)^{1/2}=\left(\int_{\Omega} u\cdot \overline{u}+\int_{\Omega} Eu: \overline{Eu}\right)^{1/2}.$$ We call $\mathcal{K}_0(\Omega)$ the closure of $C^{\infty}_0(\Omega,\mathbb{C}^N)$ with respect to the norm of $\mathcal{K}(\Omega)$.
By first Korn inequality, see for instance \cite{Ol-et-al}, we have that
$\mathcal{K}_0(\Omega)=H^1_0(\Omega,\mathbb{C}^N)$, with equivalent norms. In fact, first Korn inequality states that,
for any open set $\Omega$,
\begin{equation}\label{firstKorn}
\|\nabla u\|^2_{L^2(\Omega)}\leq 2\|E u\|^2_{L^2(\Omega)}\quad\text{for any }u\in H^1_0(\Omega,\mathbb{C}^N).
\end{equation}
By second Korn inequality, see again \cite{Ol-et-al},
provided $\Omega$ is smooth enough, for instance if $\Omega$ is a Lipschitz bounded open set, we also have that
$\mathcal{K}(\Omega)=H^1(\Omega,\mathbb{C}^N)$, with equivalent norms. Consequently, $\mathcal{K}_{loc}(\Omega)=H^1_{loc}(\Omega,\mathbb{C}^N)$.

It is well-known that, by Helmholtz decomposition, any weak solution $u$ to \eqref{Navier} can be written as
the sum of a \emph{longitudinal wave} $u_p$ and a \emph{transversal wave} $u_s$, where $u_p$ and $u_s$ are solutions
to \eqref{Navier}. Namely, we set
$$u_p=-
\frac{\nabla^T\mathrm{div}(u)}{\omega_p^2},\qquad  \omega_p^2=\frac{\rho\omega^2}{\lambda+2\mu}.$$
We note that $\chi=-\mathrm{div}(u)$ 
is a scalar weak solution to the Helmholtz equation $\Delta \chi+\omega_p^2\chi=0$ in $\Omega$, hence
$u_p$ is a vector-valued weak solution to the same Helmholtz equation
\begin{equation}\label{longwave}
\Delta u_p+\omega_p^2 u_p=0\qquad  \text{in }\Omega.
\end{equation}
If we set
$$u_s=\frac{\nabla^T\mathrm{div}(u)-\Delta u}{\omega_s^2}=
\frac{\mathrm{div}((\nabla u)^T-\nabla u)}{\omega_s^2},\qquad  \omega_s^2=\frac{\rho\omega^2}{\mu},$$
it is not difficult to show that
$u=u_p+u_s$ and
 $u_s$ 
is a vector-valued weak solution to another Helmholtz equation
\begin{equation}\label{transwave}
\Delta u_s+\omega_s^2u_s=0\qquad  \text{in }\Omega.
\end{equation}
We note that
$$u_s=\frac{\mathrm{curl}(\mathrm{curl}(u))}{\omega_s^2}\quad\text{if }N=3,\qquad
u_s=\frac{-Q\nabla^T (\mathrm{curl}_2(u))}{\omega_s^2}\quad\text{if }N=2,$$
where $Q=\left[\begin{smallmatrix} 0 & -1\\ 1 &0\end{smallmatrix}\right]$ and
$\mathrm{curl}_2(u)=u^2_1- u^1_2$ is the two-dimensional $\mathrm{curl}$ of $u$.

Since $\mathrm{div}(u_s)=0$ and $(\nabla u_p)^T-\nabla u_p=0$, we have $(u_p)_p=u_p$, $(u_p)_s=0$,
$(u_s)_s=u_s$ and $(u_s)_p=0$.

If $\Omega$ is an exterior domain, we say that $u$, a solution to \eqref{Navier}, is \emph{radiating} or {outgoing} if it satisfies the \emph{Kupradze radiation conditions}
\begin{equation}\label{kupradze}
\begin{array}{l}
\displaystyle{\lim_{r\to+\infty}r^{(N-1)/2}\left(\frac{\partial u_p}{\partial r}-\rmi \omega_p u_p\right)=0}
\\
\displaystyle{\lim_{r\to+\infty}r^{(N-1)/2}\left(\frac{\partial u_s}{\partial r}-\rmi \omega_s u_s\right)=0}\\
\end{array}\qquad r=\|x\|
\end{equation}
where the limits have to be intended as uniform in any direction. In other words, $u_p$ and $u_s$ satisfy the Sommerfeld radiation condition and, therefore, are radiating solutions to their corresponding Helmholtz equations.

For any bounded open set $\Omega$ and $u$ solution to \eqref{Navier}, the \emph{surface traction} $Tu$ is
$$Tu=\sigma(u)\nu =
2\mu Eu\nu+\lambda \mathrm{div}(u)\nu=
2\mu \nabla u\nu+\lambda\mathrm{div}(u)\nu+
\mu ((\nabla u)^T-\nabla u)\nu
\qquad\text{on }\partial \Omega,$$
$\nu$ being the exterior normal to $\Omega$, which we assume to be a column vector.
In particular, if $u$ and $\partial \Omega$ are smooth enough, we have that
$$Tu=\sigma(u)\nu =
2\mu \nabla u_p\nu
+2\mu \nabla u_s\nu
+\lambda\mathrm{div}(u_p)\nu+
\mu ((\nabla u_s)^T-\nabla u_s)\nu
\qquad\text{on }\partial \Omega.$$

For any $k>0$, let $\phi_{k}$ be the fundamental solution to the Helmholtz equation $\Delta u+k^2u =0$ which is given by
$$\phi_k(x,y)=\frac{\rmi}{4}\left(\frac{k}{2\pi\|x-y\|}\right)^{(N-2)/2}H^{(1)}_{(N-2)/2}(k\|x-y\|)
\quad\text{for any }x,\ y\in\mathbb{R}^N,\ x\neq y.$$
For any real $s \geq 0$, $H^{(1)}_{s}$ denotes the Hankel function of first kind and order $s$.
We also remark that for $N=2,3$ this reduces to the well-known formulas$$\phi_k(x,y)=\frac{\rme^{\rmi k\|x-y\|}}{4\pi\|x-y\|}\quad\text{for any }x,\ y\in\mathbb{R}^3,\ x\neq y,$$
and
$$\phi_k(x,y)=\frac{\rmi}{4}H^{(1)}_{0}(k\|x-y\|)\quad\text{for any }x,\ y\in\mathbb{R}^2,\ x\neq y.$$

Then the fundamental solution to the Navier equation is given by, for any $x,\ y\in\mathbb{R}^N$, $x\neq y$,
\begin{equation}\label{fund}
\Phi(x,y)=\frac{1}{\mu}\phi_{\omega_s}(x,y)I_N+\frac{1}{\rho\omega^2}\nabla_y\nabla_y^T\left[\phi_{\omega_s}(x,y)-\phi_{\omega_p}(x,y)\right].
\end{equation}
Here derivatives are meant in the sense of distributions over the whole $\mathbb{R}^N$ and the Navier equation is applied to $\Phi$ row by row. For $x\neq y$, we have
$$\Phi(x,y)=\frac{1}{\mu}\phi_{\omega_s}(x,y)I_N+\frac{1}{\rho\omega^2}\nabla_x\nabla_x^T\left[\phi_{\omega_s}(x,y)-\phi_{\omega_p}(x,y)\right]$$
as well. We also note that $\Phi=\Phi^T$.

For any bounded domain $\Omega$ and $u$ solution to \eqref{Navier}, provided $\Omega$ and $u$ are smooth enough, we have
for any $x\in \Omega$
\begin{equation}\label{green}
u(x)=\int_{\partial \Omega}\left(\Phi(x,y)\left[Tu(y)\right]-\left[T_y\Phi(x,y)\right] u(y)\right)\,d\sigma(y)
\end{equation}
where $T$ is applied to $\Phi$ row by row. Regarding regularity, it is enough that $\Omega$ is of class $C^2$, $u\in C^2(\Omega)\cap C(\overline{\Omega})$ such that $Tu$ exists as a uniform limit on $\partial \Omega$, namely 
$$Tu(x)=\lim_{h\to 0^+}\sigma(u)(x-h\nu(x))\nu(x)\qquad\text{for any }x\in \partial \Omega$$
where the limit is uniform with respect to $x\in \partial \Omega$ and $\nu(x)$ is the exterior normal at $x$.

If $\Omega$ is an exterior domain and $u$ is an outgoing solution to \eqref{Navier}, then we still have
for any $x\in \Omega$
\begin{equation}\label{green-unbounded}
u(x)=\int_{\partial \Omega}\left(\Phi(x,y)\left[Tu(y)\right]-\left[T_y\Phi(x,y)\right] u(y)\right)d\,\sigma(y)
\end{equation}
since the contribution at infinity is zero due to the Kupradze radiation condition satisfied by $u$ and the corresponding properties of $\Phi$.

By the well-known asymptotic properties of radiating solutions to Helmholtz equations applied to $u_p$ and to $u_s$, we infer that $u(x)=u_p(x)+u_s(x)$ satisfies
\begin{equation}\label{decay+far-field}
u(x)=
\frac{\rme^{\rmi \omega_p\|x\|}}{\|x\|^{(N-1)/2}}
U_p(\hat{x})+\frac{\rme^{\rmi \omega_s\|x\|}}{\|x\|^{(N-1)/2}}
U_s(\hat{x})
+O\left(\frac{1}{\|x\|^{(N+1)/2}}\right),
\end{equation}
as $\|x\|$ goes to $+\infty$, uniformly in all directions $\hat{x}=x/\|x\|\in \mathbb{S}^{N-1}$.
The $\mathbb{C}^N$-valued functions $U_p$ and $U_s$ are defined on $\mathbb{S}^{N-1}$ and are referred to as the \emph{longitudinal part} and the \emph{transversal part} of the
\emph{far-field pattern} $U=(U_p,U_s)$ of the field $u$, respectively.

By \eqref{green-unbounded}, the longitudinal part $U_p$ is orthogonal to $\mathbb{S}^{N-1}$, that is,
$U_p(\hat{x})=u^{\infty}_p(\hat{x})\hat{x}$ for any $\hat{x}\in \mathbb{S}^{N-1}$ for a suitable complex-valued function
$u^{\infty}_p$ defined on $\mathbb{S}^{N-1}$. On the other hand, the transversal part $U_s$ is tangential to $\mathbb{S}^{N-1}$, that is, $U_s(\hat{x})\cdot \hat{x}=0$ for any $\hat{x}\in \mathbb{S}^{N-1}$.

Therefore, if we consider the normal and tangential component, with respect to $\mathbb{S}^{N-1}$, of $u$, that is,
$$u(x)=u_N(x)+u_T(x),$$
where for any $x\in\Omega$ we have that $u_N(x)$ is proportional to $\hat{x}=x/\|x\|$ while 
$u_T(x)$ is orthogonal to $\hat{x}$, we conclude that
$$u_N(x)=\frac{\rme^{\rmi \omega_p\|x\|}}{\|x\|^{(N-1)/2}}
U_p(\hat{x})
+O\left(\frac{1}{\|x\|^{(N+1)/2}}\right),
$$
and
$$u_T(x)=\frac{\rme^{\rmi \omega_s\|x\|}}{\|x\|^{(N-1)/2}}
U_s(\hat{x})
+O\left(\frac{1}{\|x\|^{(N+1)/2}}\right),
$$
as $\|x\|$ goes to $+\infty$, uniformly in all directions. Thus, measuring the asymptotic behaviour of $u$, as 
$\|x\|$ goes to $+\infty$, corresponds to measuring both the longitudinal part and the transversal part of the
far-field pattern of $u$.

Let us send a so-called \emph{incident wave}, that is, a time-harmonic wave whose field $u^{inc}$ is an entire solution to \eqref{Navier}. Typically, the incident wave is
a \emph{plane} wave obtained by a linear combination of
a \emph{longitudinal plane wave}
\begin{equation}\label{planarlong}
u^{inc}_p(x)=d\rme^{\rmi \omega_pd\cdot x}=\frac{\nabla^T(\rme^{\rmi \omega_pd\cdot x})}{\rmi\omega_p},\qquad x\in \mathbb{R}^N,
\end{equation}
where $d\in \mathbb{S}^{N-1}$ is the \emph{direction of propagation}, and a \emph{transversal plane wave}
\begin{equation}\label{planartrans}
u^{inc}_s(x)=p\rme^{\rmi \omega_sd\cdot x},\qquad x\in \mathbb{R}^N,
\end{equation}
where $p\in \mathbb{C}^N\backslash \{0\}$ is a unitary vector orthogonal to $d$. For example, one can choose, if $N=3$,
$$u^{inc}_s(x)=c\rme^{\rmi \omega_sd\cdot x}(d_2-d_3,d_3-d_1,d_1-d_2)^T=
c\frac{\mathrm{curl}\left(\rme^{\rmi \omega_sd\cdot x}(1,1,1)^T\right)}{\rmi\omega_s},\qquad x\in \mathbb{R}^3,$$
whereas, if $N=2$,
$$u^{inc}_s(x)=-cQd\rme^{\rmi \omega_sd\cdot x}=c\frac{-Q\nabla^T(\rme^{\rmi \omega_sd\cdot x})}{\rmi\omega_s},\qquad x\in \mathbb{R}^2,$$
where $c\in\mathbb{C}\backslash\{0\}$ is a suitable constant.

Namely, we consider
\begin{equation}\label{planarcomb}
u^{inc}(x)=c_pu^{inc}_p(x)+c_su^{inc}_s(x),\qquad x\in \mathbb{R}^N
\end{equation}
for some $d\in \mathbb{S}^{N-1}$, $p\in \mathbb{C}^{N}\backslash \{0\}$ such that $\|p\|=1$ and $p$ is orthogonal to $d$, and $c_p,c_s\in\mathbb{C}$ such that $|c_p|^2+|c_s|^2=1$, in such a way to have
\begin{equation}\label{planarcomb-prop}
\|u^{inc}(x)\|=1,\qquad x\in \mathbb{R}^N
\end{equation}

The presence of an impenetrable object, that is, of a scatterer $K$, inside the medium perturbs the incident wave by creating the \emph{scattered} or \emph{reflected wave}, given by the field $u^{scat}$. The total wave is the superposition of the incident and the scattered waves and its field is denoted by $u$. The total field $u$ solves \eqref{Navier} in $\Omega=\mathbb{R}^N\backslash K$ and satisfies a boundary condition on $\partial K$ that depends on the nature of the scatterer, namely, if $K$ is a so-called \emph{rigid scatterer}, a homogeneous Dirichlet boundary condition 
\begin{equation}\label{BC1}
u=0\qquad\text{on }\partial K,
\end{equation}
or, if $K$ is a so-called \emph{cavity}, a homogeneous Neumann boundary condition
\begin{equation}\label{BC2}
Tu=\sigma(u)\nu =0\qquad\text{on }\partial K, 
\end{equation}
$\nu$ being the exterior normal to $\Omega$.

Finally, being $\Omega=\mathbb{R}^N\backslash K$ unbounded, a condition at infinity has to be imposed. We require the scattered wave to be outgoing.
Summarising, the total field $u$ solves the following exterior boundary value problem
\begin{equation}\label{scatteringpbmrigid}
\left\{
\begin{array}{ll}
u=u^{inc}+u^{scat} &\text{in }\Omega=\mathbb{R}^N\backslash K\\
\mu\Delta u+(\lambda+\mu)\nabla(\mathrm{div} (u))+\rho \omega^2 u=0 &\text{in }\Omega\\
u=0 &\text{on }\partial \Omega=\partial K\\
u^{scat}\text{ satisfies }\eqref{kupradze}
\end{array}
\right.
\end{equation}
if $K$ is a rigid scatterer, and 
\begin{equation}\label{scatteringpbmcavity}
\left\{
\begin{array}{ll}
u=u^{inc}+u^{scat} &\text{in }\Omega\\
\mu\Delta u+(\lambda+\mu)\nabla(\mathrm{div} (u))+\rho \omega^2 u=0 &\text{in }\Omega\\
Tu=0 &\text{on }\partial \Omega\\
u^{scat}\text{ satisfies }\eqref{kupradze}
\end{array}
\right.
\end{equation}
if $K$ is a cavity.

The weak formulation of \eqref{scatteringpbmrigid} is the following. Assume that $K\subset B_R(0)$ for some $R>0$. Then we look for $u$ belonging to $H^1(B_r(0)\backslash K,\mathbb{C}^N)$ for any $r> R$ such that $u=u^{inc}+u^{scat}$ solves \eqref{Navier} in the weak sense and $u^{scat}$ satisfies the condition at infinity given by \eqref{kupradze}. Finally, for what concerns the boundary condition \eqref{BC1} on $\partial K$, we require that $u=0$ on $\partial K$ in a weak sense, that is,
$\chi u\in H^1_0(B_r(0)\backslash K,\mathbb{C}^N)$ for any $r>R$ and any $\chi\in C^{\infty}_0(B_r(0),\mathbb{R})$ such that $\chi=1$ on $B_R(0)$.

The weak formulation of \eqref{scatteringpbmcavity} is the following. Assume that $K\subset B_R(0)$ for some $R>0$. Then we look for $u$ belonging to $\mathcal{K}(B_r(0)\backslash K,\mathbb{C}^N)$ for any $r> R$ such that $u=u^{inc}+u^{scat}$ solves \eqref{Navier} in the weak sense and $u^{scat}$ satisfies the condition at infinity given by \eqref{kupradze}. Finally, about the boundary condition \eqref{BC2} on $\partial K$, we require that
for any $r>R$ and any $v\in \mathcal{K}(B_r(0)\backslash K,\mathbb{C}^N)$ with compact support contained in $B_r(0)$
we have
\begin{equation}\label{Navierweak2}
2\mu\int_{\Omega}Eu : \overline{Ev}+\lambda\int_{\Omega}\mathrm{div} (u) \overline{\mathrm{div} (v)} - \rho\omega^2\int_{\Omega}u\cdot\overline{v}=0.
\end{equation}

In both cases, if $u^{inc}=0$, then we have that
for any $r>R$
\begin{equation}\label{byparts}
\Im\left(\int_{\partial B_r(0)}(Tu)\cdot\overline{u}\right)=0,
\end{equation}
hence, by the asymptotic behaviour of outgoing solutions, we infer that $u=u^{scat}=0$ in $\Omega$. In other words,
\eqref{scatteringpbmrigid} and
 \eqref{scatteringpbmcavity} admit at most one solution, so uniqueness follows. Concerning existence, this can be established by layer potential techniques provided $K$, or $\Omega$, is regular enough, say of class $C^2$. In fact,
for any exterior domain $\Omega$ of class $C^2$ and for any $\varphi\in C(\partial \Omega,\mathbb{C}^3)$ we define the single- and double-layer potentials with density $\varphi$ at any $x\in\Omega$ as follows
\begin{equation}\label{single-double}
\mathcal{S}(\varphi)(x)=
\int_{\partial \Omega}\Phi(x,y)\varphi(y)\,d\sigma(y)
\quad\text{and}\quad 
\mathcal{D}(\varphi)(x)=
\int_{\partial \Omega}\left[T_y\Phi(x,y)\right] \varphi(y)\,d\sigma(y).
\end{equation}
We observe that $\mathcal{S}(\varphi)$ and $\mathcal{D}(\varphi)$ are outgoing solutions to \eqref{Navier} in $\Omega$.
By carefully exploiting the properties of the potentials on $\partial\Omega$, the following existence (and uniqueness) result can be proved.

\begin{teo}\label{existenceteo}
Assume that $K$ is an obstacle such that $D=\overset{\circ}{K}$ is an open set of class $C^2$. Then, for any $u^{inc}$ entire solution to \eqref{Navier}, \eqref{scatteringpbmrigid} and
 \eqref{scatteringpbmcavity} admit one solution.
\end{teo}

\proof{.} This is a classical result, see for instance \cite{Kup2}. Actually, the regularity of $D$ can be relaxed up to Lipschitz, see for instance \cite[Corollary~2,3]{Bao-et-al}.\cvd

\smallskip

We shall need the following regularity result, whose proof is postponed to the appendix.

\begin{teo}\label{regularityteo}
Let us fix positive constants $r$, $L$, $R$ and $\alpha$, $0<\alpha<1$. Let us also fix the coefficients $\mu>0$, $\lambda$ such that
$2\mu+\lambda>0$, $\rho>0$ and $\omega>0$.
Assume that $K\subset B_R(0)$ is an obstacle such that $D=\overset{\circ}{K}$ is an open set of class $C^{2,\alpha}$ with constants $r$ and $L$. Let $u^{inc}$ be as in \eqref{planarcomb} such that \eqref{planarcomb-prop} is satisfied.

Let $u$ be the solution to \eqref{scatteringpbmrigid}. Then there exists a constant $\tilde{C}_0$, depending on $r$, $L$, $R$, $\alpha$ and the coefficients only, such that
\begin{equation}\label{regularity}
\|u\|_{C^{2}(\overline{\Omega})}\leq \tilde{C}_0.
\end{equation}

Moreover, there exists a constant $\tilde{C}_1$, depending on $r$, $L$, $R$, $\alpha$ and the coefficients only, such that
\begin{equation}\label{decay}
\left\|u^{scat}(x)\right\|\leq  \left\|u^{scat}_s(x)\right\|+\left\|u^{scat}_p(x)\right\|\leq \frac{\tilde{C}_1}{\|x\|^{(N-1)/2}}\qquad \text{for any }\|x\|\geq R+1.
\end{equation}
\end{teo}

We conclude this part with the 
following regularity result and a three-spheres inequality for the Helmholtz equation.

\begin{lem}\label{regullemma}
Let $u$ be a solution to \eqref{Navier} in $B_s$, with $0<s\leq s_0$. Then there exists a constant $D_0$, depending on the coefficients of \eqref{Navier} and on $s_0$ only,    such that
\begin{equation}\label{regestimate-a}
\|u_p\|_{L^{\infty}(B_{s/8})},\ \|u_s\|_{L^{\infty}(B_{s/8})}\leq \frac{D_0}{s^{(N+4)/2}} \|u\|_{L^2(B_s)}
\end{equation}
and
\begin{equation}\label{regestimate-b}
\|u_p\|_{L^{\infty}(B_{s/4})},\ \|u_s\|_{L^{\infty}(B_{s/4})}\leq \frac{D_0}{s^{(N+2)/2}} \|\nabla u\|_{L^2(B_s)}.
\end{equation}
\end{lem}

\proof{.} First of all we use a Caccioppoli inequality to estimate the $H^1$-norm of $u$ in a smaller ball.

Namely, let $\chi\in C^{\infty}_0(B_s)$ be such that $0\leq \chi\leq 1$ everywhere and $\chi=1$ on $B_{3s/4}$. We can assume that $\|\nabla \chi\|\leq C/s$ everywhere for some absolute constant $C$.
Then we apply the weak formulation of \eqref{Navier} to $v=\chi^2 u$ and obtain
$$2\mu\int_{B_s}Eu : \overline{Ev}+\lambda\int_{B_s}\mathrm{div} (u) \overline{\mathrm{div} (v)} =
\rho\omega^2\int_{B_s}\chi^2\|u\|^2.$$
But 
\begin{multline*}
2\mu\int_{B_s}Eu : \overline{Ev}+\lambda\int_{B_s}\mathrm{div} (u) \overline{\mathrm{div} (v)}=
2\mu\int_{B_s}\chi^2 \|Eu\|^2+\lambda\int_{B_s}\chi^2|\mathrm{div} (u)|^2\\+
\left(2\mu\int_{B_s}\chi Eu:(\chi_j\overline{u}^i+\chi_i\overline{u}^j)+2\lambda\int_{B_s}\chi \mathrm{div} (u)
(\chi_j\overline{u}^j)\right),
\end{multline*}
hence, by Cauchy inequality, we infer that
\begin{multline*}
2\mu\int_{B_s}\chi^2 \|Eu\|^2+\lambda\int_{B_s}\chi^2|\mathrm{div} (u)|^2\\
\leq
\mu\int_{B_s}\chi^2 \|Eu\|^2+\frac{\lambda}{2}\int_{B_s}\chi^2|\mathrm{div} (u)|^2
+\left(\frac{C^2}{s^2}C_1+\rho\omega^2\right)\int_{B_s} \|u\|^2
\end{multline*}
for a constant $C_1$ depending on $\mu$ and $\lambda$ only. We conclude that
\begin{multline*}
\frac{\min\{2\mu,2\mu+\lambda\}}{2}\int_{B_{3s/4}}\|Eu\|^2\leq 
\mu\int_{B_{3s/4}} \|Eu\|^2+\frac{\lambda}{2}\int_{B_{3s/4}}|\mathrm{div} (u)|^2\\\leq
\mu\int_{B_s}\chi^2 \|Eu\|^2+\frac{\lambda}{2}\int_{B_s}\chi^2|\mathrm{div} (u)|^2\leq
\left(\frac{C^2}{s^2}C_1+\rho\omega^2\right)\int_{B_s} \|u\|^2.
\end{multline*}

Let $\tilde{\chi}\in C^{\infty}_0(B_{3s/4})$ be such that $0\leq \tilde{\chi}\leq 1$ everywhere and $\tilde{\chi}=1$ on $B_{s/2}$. We can assume that $\|\nabla \tilde{\chi}\|\leq C/s$ everywhere for some absolute constant $C$. Then $\tilde{\chi} u\in\mathcal{K}_0(B_{3s/4})=
H^1_0(B_{3s/4})$, so by \eqref{firstKorn},
\begin{multline*}
\int_{B_{s/2}}\|\nabla u\|^2\leq \int_{B_{3s/4}}\|\nabla (\tilde{\chi} u)\|^2\leq
2\int_{B_{3s/4}}\|E (\tilde{\chi} u)\|^2\\=2\int_{B_{3s/4}}\left\|\tilde{\chi} Eu + \frac12(\tilde{\chi}_j u^i+\tilde{\chi}_i u^j)\right\|^2\leq4\left(\int_{B_{3s/4}}\|Eu\|^2+\frac{C^2}{s^2}C_2\int_{B_{3s/4}}\|u\|^2
\right)
\end{multline*}
where $C_2$ is another absolute constant. We conclude that
\begin{equation}\label{stepone}
\|\nabla u\|_{L^2(B_{s/2})}\leq \frac{C_3}{s}\|u\|_{L^2(B_s)}
\end{equation}
for a constant $C_3$ depending on the coefficients and on $s_0$ only.

Since for any $j\in\{1,\ldots,N\}$, $u_j$ still solve \eqref{Navier}, we can repeat the procedure above and prove that
\begin{equation}\label{steptwo}
\|D^2 u\|_{L^2(B_{s/4})}\leq \frac{C_4}{s^2}\|u\|_{L^2(B_s)},
\end{equation}
consequently
\begin{equation}\label{steptwobis}
\|u_p\|_{L^2(B_{s/4})}, \|u_s\|_{L^2(B_{s/4})}\leq \frac{C_5}{s^2}\|u\|_{L^2(B_s)},
\end{equation}
with $C_4$ and $C_5$ still depending on the coefficients and on $s_0$ only.

The last step is to estimate the $L^{\infty}$-norm by the $L^2$-norm for a solution to a Helmholtz equation. This is a standard estimate, see for instance \cite[Theorem~8.17]{G-T}, since we have
\begin{equation}\label{infty-2}
\|u_p\|_{L^{\infty}(B_{s/8})}, \|u_s\|_{L^{\infty}(B_{s/8})}\leq \frac{C_6}{s^{N/2}}\|u_p\|_{L^2(B_{s/4})}, \|u_s\|_{L^2(B_{s/4})},
\end{equation}
respectively, with $C_6$ depending on the coefficients and on $s_0$ only. The proof can now be easily concluded.\cvd

\smallskip

\begin{lem}\label{3sphereslemma}
There exist positive  constants $\tilde{s}_0$, $\tilde{C}$ and $\tilde{c}_1$, $0<\tilde{c}_1<1$,
depending on $k$ only, such that for every 
$0<s_1<s<s_2\leq \tilde{s}_0$
and any function $u$ such that
$$\Delta u+k^2u=0\quad\text{in }B_{s_2},$$
we have, for any $t$, $s<t<s_2$,
\begin{equation}\label{3spheres}
\|u\|_{L^{\infty}(B_{s})}\leq \tilde{C}(1-(s/t))^{-N/2}
\|u\|^{1-\beta}_{L^{\infty}(B_{s_2})}\|u\|^{\beta}_{L^{\infty}(B_{s_1})},\end{equation}
for some $\beta$ such that
\begin{equation}\label{3spherescoeff}
\tilde{c}_1\left(\log(s_2/t)\right)\big/\left(\log(s_2/s_1)\right)
\leq\beta\leq 1-\tilde{c}_1\left(\log(t/s_1)\right)\big/\left(\log(s_2/s_1)\right).
\end{equation}
\end{lem}

\proof{.} It follows by the results of \cite{Bru}.\cvd

\smallskip

\subsection{Friedrichs inequality}\label{Friedrichs inequality}

Let $\Omega\subset \mathbb{R}^N$ be an open and bounded set.  Under suitable assumptions on $\Omega$ and $u$,
 a function defined on $\overline{\Omega}$, the inequality proved by Friedrichs, \cite{Fri}, is
\begin{equation}\label{Friinitial}
\|u\|_{L^2(\Omega)}\leq C\left[\|\nabla u\|_{L^2(\Omega)}+\|u\|_{L^2(\partial\Omega)}\right].
\end{equation}
with a constant $C$ not depending on $u$. This estimate is actually a straightforward consequence of a much more general estimate proved by Maz'ya, \cite{M0}, which is the following
\begin{equation}\label{fri0}
\|u\|_{L^{N/(N-1)}(\Omega)}\leq C(N)\left[\|\nabla u\|_{L^1(\Omega)}+\|u\|_{L^1(\partial\Omega)}\right]\quad \text{for any }u\in C(\overline{\Omega})\cap W^{1,1}(\Omega).
\end{equation}
The importance of this estimate is that it holds independently of the regularity of $\Omega$ and that the constant $C(N)$ is optimal and depends on $N$ only, it is actually the one of the isoperimetric inequality, that is,
\begin{equation}\label{constant}
C(N)=\frac{|B_1|^{(N-1)/N}}{\mathcal{H}^{N-1}(\partial B_1)}.
\end{equation}
For a proof of \eqref{fri0} we refer to \cite[Corollary on page~319]{M}. Actually, the Maz'ya inequality can even be generalised to functions of bounded variation, see \cite{R} for an extremely general version in this direction.

Here we just point out that \eqref{fri0} implies the classical Friedrichs inequality, which we state in the next theorem.

\begin{teo}\label{Friedrichsthm}
Let $\Omega\subset \mathbb{R}^N$ be open and bounded. Let
$$p=\frac{2N}{N+1}\qquad\text{and}\qquad q=\frac{2N}{N-1}.$$
Then
\begin{equation}\label{Friedrichs0}
\|u\|_{L^q(\Omega)}\leq 4C(N)\left[\|\nabla u\|_{L^p(\Omega)}+\|u\|_{L^2(\partial\Omega)}\right]
\quad \text{for any }u\in C(\overline{\Omega})\cap W^{1,p}(\Omega),
\end{equation}
so that 
for any $u\in C(\overline{\Omega})\cap W^{1,2}(\Omega)$
\begin{equation}\label{Friedrichs1}
\|u\|_{L^2(\Omega)}\leq 4C(N)|\Omega|^{\frac{1}{2N}}\left[|\Omega|^{\frac{1}{2N}}\|\nabla u\|_{L^2(\Omega)}+\|u\|_{L^2(\partial\Omega)}\right].
\end{equation}
\end{teo}

\proof{.} Inequality \eqref{Friedrichs0} follows by applying \eqref{fri0} to $u^2$, see the proof of \cite[Corollary~2.4]{R}
for details, whereas \eqref{Friedrichs1} is an immediate consequence of \eqref{Friedrichs0} and H\"older inequality.\cvd

\smallskip

When $u$ is vector-valued, that is,
$u\in C(\overline{\Omega},\mathbb{C}^N)\cap W^{1,2}(\Omega,\mathbb{C}^N)$, we still have
\begin{equation}\label{Friedrichs}
\|u\|_{L^2(\Omega)}\leq 4C(N)|\Omega|^{\frac{1}{2N}}\left[|\Omega|^{\frac{1}{2N}}\|\nabla u\|_{L^2(\Omega)}+\|u\|_{L^2(\partial\Omega)}\right].
\end{equation}

\section{The main result}\label{mainresultsec}

We begin by setting the hypotheses. Let us fix constants $r>0$, $L>0$, $R>0$,
 $\alpha$ with $0<\alpha<1$,
 $\mu>0$, $\lambda$ with $2\mu+\lambda>0$, $\rho>0$, $\omega>0$. Finally, we fix $H_0$ such that
\begin{equation}\label{closenessconstant}
0<H_0<H_1=\left(\frac{\min\{2\mu,2\mu+\lambda\}}{64 C(N)^2\rho\omega^2}\right)^{N/2}=
\left( \frac{\min\{\omega_p^{-1},\sqrt{2}\omega_s^{-1}\}}{8C(N)}\right)^N  
\end{equation}
where $C(N)$ is the absolute constant appearing in \eqref{constant}.
We refer to these constants as the \emph{a priori data}.

First of all, we fix $D^+$, an open set which is Lipschitz with constants $r$ and $L$ and such that $K^+=\overline{D^+}$ is an obstacle contained in $B_R(0)$. We call $\Omega^+=\mathbb{R}^N\backslash K^+$.

Let $D$ and $D'$ be two open sets of class $C^{2,\alpha}$
with constants $r$ and $L$ such that
$K=\overline{D}$ and $K'=\overline{D'}$ are obstacles contained in $D^+$.
We call $\Omega=\mathbb{R}^N\backslash K$
and $\Omega'=\mathbb{R}^N\backslash K'$, the corresponding exterior domains. We also use the notation $\Omega^{ext}$ to denote the unbounded connected component of $\mathbb{R}^N\backslash (K\cup K')$ and call $\Gamma=\partial \Omega^{ext}$ and $\Omega^{int}=\mathbb{R}^N\backslash\overline{\Omega^{ext}}$.
We clearly have $\Omega^+\subset\Omega^{ext}$ and $\Omega^{int}\subset D^+$.

Let $u^{inc}$ be as in \eqref{planarcomb} such that \eqref{planarcomb-prop} is satisfied. Let $u$ and $u^{scat}$ be the solution to \eqref{scatteringpbmrigid} and let $u'$ and $(u')^{scat}$ be the solution to \eqref{scatteringpbmrigid} with $K$ replaced by $K'$.

Let $U=(U_p,U_s)$ be the far-field pattern of $u^{scat}$ and $U'=(U'_p,U'_s)$ be the far-field pattern of $(u')^{scat}$, respectively.

We measure the difference between two obstacles $K$ and $K'$ by using the \emph{Hausdorff distance} $d_H$ which is given by
$$d_H(K,K')=\max\left\{\sup_{x\in K}\mathrm{dist}(x,K'),\sup_{x\in K'}\mathrm{dist}(x,K)\right\}.$$

Then we have the following stability result.

\begin{teo}\label{mainthm}
Under the previous notation and assumptions, we further assume that the following \emph{closeness condition} holds
\begin{equation}\label{closeness}
\left|D^+\backslash (K\cap K')\right|\leq H_0.
\end{equation}

Then there exist positive constants $\hat{\varepsilon}_0$, $0<\hat{\varepsilon}_0\leq \rme^{-\rme}/2$, $\hat{C}$ and $\beta$, depending on the a priori data only, such that for any $\varepsilon_0$, $0<\varepsilon_0\leq \hat{\varepsilon}_0$,
if
\begin{equation}\label{farfield-error}
\|U-U'\|_{L^2(\mathbb{S}^{N-1},\mathbb{C}^N\times\mathbb{C}^N)}\leq \varepsilon_0,
\end{equation}
then
\begin{equation}
d_H(K,K')\leq \hat{C}\left(\log\left(\log(1/\varepsilon_0)\right)\right)^{-\beta}.
\end{equation}
\end{teo}

\begin{oss} As it will appear clear in the proof, about the incident wave, we just need conditions that allow \eqref{farawaypoint} and \eqref{farawaypoint2} to be satisfied.
Therefore, other suitable incident waves may be used. For instance, another common choice is to use a point source wave. However, in such a case, one needs to consider a point source which is far enough from the unknown obstacle and choose $x_3$ in \eqref{farawaypoint} and \eqref{farawaypoint2} relatively close to the point source. The analysis would therefore require other technicalities that we decided not to tackle in this paper.
\end{oss}

\section{Proof of Theorem~\ref{mainthm}.}\label{proofsec}

Let $K$ and $K'$ be any two scatterers satisfying the hypotheses of Theorem~\ref{mainthm}. We state a few properties of $K$, as well as of $K'$.
First of all we note that the number of connected components of $K$ is bounded by a constant depending on $r$, $L$ and $R$ only. We also have that $\mathcal{H}^{N-1}(\partial K)$ is bounded by a constant depending on $r$, $L$ and $R$ only.
Moreover, there exists a constant $C_1$, depending on $r$, $L$ and $R$ only, such that for any $h$, $0\leq h\leq 1$, we have
\begin{equation}\label{boundary33}
|\overline{B_{h}(\partial K)}|\leq C_1 h.
\end{equation}

By \cite[Corollary~2.3 and Proposition~2.1]{LPRX}, there exist two positive constants $c_1$ and $t_1$, depending on $r$, $L$ and $R$ only, such that the following holds.
For any $t>0$,
if $x_1$, $x_2\in\mathbb{R}^N$ are such that
$B_t(x_1)$ and $B_t(x_2)$ are contained in $\mathbb{R}^N\backslash K$,
then we can find a smooth (for instance piecewise $C^1$) curve $\gamma$ connecting
$x_1$ to $x_2$ so that $\overline{B_{\delta(t)}(\gamma)}$ is contained in$\mathbb{R}^N\backslash K$ as well, where
\begin{equation}\label{deltaconnect}
\delta(t)=\min\{c_1 t, t_1\}\quad\text{for any }t>0.
\end{equation}

We measure the distance between $K$ and $K'$ by
\begin{equation}\label{ddefin}
d=\max\left\{\sup_{x\in\partial K\backslash K'}\mathrm{dist}(x,\partial K'),
\sup_{x\in\partial K'\backslash K}\mathrm{dist}(x,\partial K)\right\}
\end{equation}
or
\begin{equation}\label{dddefin}
\hat{d}=d_H(\partial K,\partial K')\quad\text{or}\quad\tilde{d}=d_H(K,K').
\end{equation}
We obviously have $d,\ \hat{d},\ \tilde{d}\leq 2R$. 
The relationship between these quantities is investigated in detail in \cite[Section~2]{LPRX} under much more general conditions. Here we just use that, in particular by \cite[Corollary~2.3 and Proposition~2.1]{LPRX}, we have
\begin{equation}\label{distancesoppositebis}
C_2 d\leq C_2\hat{d}\leq \tilde{d}\leq C_3d\leq C_3\hat{d},
\end{equation}
where $C_2$ and $C_3$ are positive constants depending on $r$, $L$ and $R$ only.

Let us note that all the above properties are valid even if we assume that $D$ and $D'$ are just Lipschitz with constants $r$ and $L$. In particular \eqref{boundary33} holds with $K$ replaced by $K^+$ as well.

If $D$ is Lipschitz with constants $r$ and $L$, $D$ and $\Omega$ satisfy a \emph{uniform interior cone property}, that is, there exist constants $r_0>0$ and $\theta_0$, $0<\theta_0<\pi/2$, depending on $r$ and $L$ only, such that for any $x\in \partial K$ we can find a unit vector $v$
such that $\mathcal{C}(x,v,r_0,\theta_0)\subset D$ and $\mathcal{C}(x,-v,r_0,\theta_0)\subset \Omega$.
Let us also note that $v$ can be chosen constant for any $y\in\partial K$ in a neighbourhood of $x\in \partial K$ depending on $r$ and $L$ only.

Another important property, for which $D$ of class $C^{1,1}$ with constants $r$ and $L$ would be enough, is the following, see \cite[Theorem~5.7]{Del-Zol}. There exist positive constants $h_0$, $r_1$ and $L_1$, depending on $r$ and $L$ only, such that
for any $h$, $0< h\leq h_0$, the set
$$D_h=\{x\in\mathbb{R}^N:\ \mathrm{dist}(x,K)<h\}$$
is an open set of class $C^{1,1}$ with constants $r_1$ and $L_1$. Moreover,
$$\partial (D_h)=\{x\in\mathbb{R}^N:\ \mathrm{dist}(x,K)=h\}.$$
We can conclude that there exists a constant $C_4$,
depending on $r$, $L$ and $R$ only, such that for any $h$, $0\leq h\leq h_0$,
\begin{equation}\label{measureboundary}
\mathcal{H}^{N-1}(\partial (D_h))\leq C_4,
\end{equation}
where we identify $D_0$ with $D$ and $\partial (D_0)$ with $\partial D$.

The final property we need about the obstacles $K$ and $K'$ is the following.

\begin{lem}\label{RLGlemma}
Assume that $D$ and $D'$ are $C^{1,1}$ with constants $r$ and $L$. Then
there exists a constant $\tilde{d}_0$, depending on $r$, $L$ and $R$ only, such that if
$$\tilde{d}=d_H(K,K')\leq \tilde{d}_0,$$
then $\Omega^{ext}$ satisfies a uniform interior cone property, with constants $\tilde{r}_0$ and $\tilde{\theta}_0$ depending on $r$ and $L$ only.\end{lem}

\proof{.} It immediately follows from \cite[Lemma~8.1]{Al-etal}. We just note that $D$ and $D'$ belonging to $C^{1,\alpha}$, with $0<\alpha< 1$, with constants $r$ and $L$ would be enough, but in this case the constants would depend on $\alpha$ as well.\cvd

\smallskip

By Theorem~\ref{regularityteo}, we have that
\begin{equation}\label{globalbound}
\|u_p(x)\|+\|u'_p(x)\|+\|u_s(x)\|+\|u'_s(x)\|\leq E\quad\text{for any }x\in \overline{\Omega^{ext}},
\end{equation}
where $E$ depends on the a priori data only it is assumed to be greater than or equal to $1$.

Finally, we fix positive $R_1$ and $\tilde{s}$ such that
$R+1+\tilde{s}\leq R_1$.
Let us fix a point $x_0\in\mathbb{R}^N$ such that
$R+1+\tilde{s}\leq \|x_0\|\leq R_1$.
For a fixed $\varepsilon$, $0<\varepsilon\leq E$, let
\begin{equation}\label{errornear}
\|u-u'\|_{L^{\infty}(B_{\tilde{s}}(x_0),\mathbb{C}^N)}\leq \varepsilon.
\end{equation}
We call $\varepsilon$ the \emph{near-field error with limited aperture}.
Let $\varepsilon_1$, $0<\varepsilon_1\leq E$, be such that 
\begin{equation}\label{errornear2}
\|u-u'\|_{L^{\infty}(B_{\|x_0\|+\tilde{s}}\backslash \overline{B_{\|x_0\|-\tilde{s}}},\mathbb{C}^N)}\leq \varepsilon_1.
\end{equation}
We call $\varepsilon_1$ the \emph{near-field error}.
Finally,
if
\begin{equation}\label{errorfar}
\|U-U'\|_{L^2(\mathbb{S}^{N-1},\mathbb{C}^N\times\mathbb{C}^N)}\leq \varepsilon_0,
\end{equation}
$\varepsilon_0$ will be referred to as the \emph{far-field error}.

By Theorem~\ref{regularityteo}, through \eqref{globalbound}, Lemma~\ref{regullemma} and an iterated application of the three-spheres inequality of Lemma~\ref{3sphereslemma} to $u_p-u'_p$ and $u_s-u'_s$,
we can find positive constants $C_5$ and $\tilde{\beta}$, $0<\tilde{\beta}<1$, depending on $E$, $R$, $\tilde{s}$, $R_1$ and the coefficients of \eqref{Navier} only, such that
\begin{equation}\label{ballannulus}
\varepsilon\leq \varepsilon_1\leq C_5\varepsilon^{\tilde{\beta}}.
\end{equation}

Moreover, 
there exist positive constants $\tilde{\varepsilon}_0\leq 1/(2\rme)$ and $C_6$, depending on 
$E$, $R$, $\tilde{s}$, $R_1$ and the coefficients of \eqref{Navier} only,such that if $0<\varepsilon_0\leq\tilde{\varepsilon}_0$ then
\begin{multline}\label{fartonear}
\|u-u'\|_{L^{\infty}(B_{\|x_0\|+\tilde{s}}\backslash \overline{B_{\|x_0\|-\tilde{s}}},\mathbb{C}^N)}\\\leq
\|u_p-u'_p\|_{L^{\infty}(B_{\|x_0\|+\tilde{s}}\backslash \overline{B_{\|x_0\|-\tilde{s}}},\mathbb{C}^N)}+
\|u_s-u'_s\|_{L^{\infty}(B_{\|x_0\|+\tilde{s}}\backslash \overline{B_{\|x_0\|-\tilde{s}}},\mathbb{C}^N)}\\
\leq
\tilde{\eta}(\varepsilon_0)=\exp\left(-C_6(-\log\varepsilon_0)^{1/2}\right),
\end{multline}
that is, possibly slightly changing $\tilde{\varepsilon}_0$
\begin{equation}\label{fartonear2}
\varepsilon\leq \varepsilon_1\leq \tilde{\eta}(\varepsilon_0)=\exp\left(-C_6(-\log\varepsilon_0)^{1/2}\right)\leq
\exp\left(-(\log(1/\varepsilon_0))^{1/4}\right).
\end{equation}
This is a classical far-field to near-field estimate, which has been first introduced in \cite{Isak92} for $N=3$, with a slight improvement in \cite{Bus}, and that can be generalised to any $N\geq 2$, see for instance Theorem~4.1 in \cite{R-S}.

We estimate the Hausdorff distance of $K$ and $K'$ in terms of $\varepsilon$. In this case, 
we need to add $R_1$ and $\tilde{s}$ to the a priori data. By \eqref{ballannulus}, the estimate in terms of $\varepsilon_1$ is clearly the same.
The estimate in terms of the far-field error $\varepsilon_0$ can be easily obtained by using \eqref{fartonear2}, noting that in this case $R_1$ and $\tilde{s}$ can be chosen as depending on $R$ only.

\label{V_sdefin}
For any $s>0$ let us call $V_{s}$ the set of points $x\in \Omega^{ext}$ such that there exists a smooth, that is, piecewise $C^1$, curve $\gamma$ connecting $x_0$ to $x$ such that $\overline{B_{s}(\gamma)}\subset \Omega^{ext}$ (see also \cite{AS} for a related argument developed in order to circumvent the case in which a  domain of interest is not reachable by
a chain of balls).
It follows that $V_{s}$ is an open subset of $\Omega^{ext}$ and we call $\Gamma_{s}$ its boundary and
$W_{s}=\mathbb{R}^N\backslash \overline{V_{s}}$. To keep the same notation, we identify $\Omega^{ext}$ with $V_0$, $\Gamma$ with $\Gamma_0$ and $\Omega^{int}$ with $W_0$. For any $0\leq s_1\leq s_2$ we clearly have
$V_{s_2}\subset V_{s_1}$ and $W_{s_1}\subset W_{s_2}$.

An important property of $\Gamma_s$ is that
\begin{equation}\label{Gammaprop1}
\Gamma_{s}\subset  \partial (D_{s})\cup\partial (D'_{s}).
\end{equation}
Moreover, by \eqref{deltaconnect} applied to $K^+$, for any $x\in \Omega^+$ whose distance from $\partial D^+$ is greater than or equal to $t>0$, we have
$x\in V_{\delta(t)}$.

We can find $s_0$, $0<s_0\leq\tilde{s}/8$, depending on the a priori data only, such that the following holds. It is smaller than or equal to $\tilde{s}_0$ in Lemma~\ref{3sphereslemma} for $k$ equal to $\omega_p$ and to $\omega_s$.
It is smaller than or equal to $h_0$. Finally, we require
that
\begin{equation}\label{boundary44}
\left|W_{s}\backslash K^+\right|\leq \frac{H_1-H_0}{2}\quad\text{for any }0\leq s\leq s_0.
\end{equation}
For this last property we use \eqref{deltaconnect} and \eqref{boundary33} applied to $K^+$. 
By \eqref{boundary44}, via \eqref{closeness}, we infer that
\begin{equation}\label{boundary55}
\left|W_{s}\backslash (K\cap K')\right|\leq \frac{H_0+H_1}{2}=\tilde{H}_0<H_1\quad\text{for any }0\leq s\leq s_0.
\end{equation}

Up to swapping $K$ with $K'$, let $x_1\in \partial K'\backslash K$ be such that
$\mathrm{dist}(x_1,\partial K)=\mathrm{dist}(x_1,K)=d$. If $x_1$ does not belong to $\Gamma$,
we can find a smooth curve $\gamma$ connecting $x_1$ with $x_0$ such that$\overline{B_{\delta(d)}(\gamma)}\subset \Omega$. But $\gamma$ needs to intersect $\partial K'\cap \Gamma$ in a point. Therefore, for a positive constant $c_2$, $0<c_2<1$, depending on $r$, $L$ and $R$ only, we can assume, without loss of generality, that there exists
$x_1\in (\partial K'\cap \Gamma)\backslash K$ such that
$\mathrm{dist}(x_1,\partial K)=\mathrm{dist}(x_1,K)\geq c_2d$.

We call $A_0$ the connected component of $\Omega^{int}\backslash K$ such that $x_1\in \partial A_0$.
For any $s$, $0<s\leq s_0$, we call
$A_{s}$ the connected component of $W_{s}\backslash K$ containing $A_0$.
For any $s$, $0\leq s\leq s_0$,
the domain $A_{s}$ satisfies the following properties. 
By \eqref{Gammaprop1},
$\partial A_{s}\subset \partial K\cup \partial (D_{s})\cup\partial (D'_{s})$. Therefore, by \eqref{measureboundary} and by \eqref{boundary55}, we have
\begin{equation}\label{measurebdry+closeness}
\mathcal{H}^{N-1}(\partial A_{s})\leq 3C_4\qquad\text{and}\qquad |A_{s}|\leq \frac{H_0+H_1}{2}=\tilde{H}_0.
\end{equation}
Moreover, by the regularity of $D$ and $D'$, we infer that there exist a point $x_2$ and a positive constant $c_3$, $0<c_3<1$ depending on $r$, $L$ and $R$ only, such that
\begin{equation}\label{startingball}
B_{c_3 d}(x_2)\subset A_0.
\end{equation}

By \eqref{planarcomb-prop} and \eqref{decay}, we can find a constant $R_2\geq R+2$, depending on the a priori data only,  and a point $x_3$ such that
\begin{equation}\label{farawaypoint}
R+2\leq\|x_3\|\leq R_2
\end{equation}
and
 \begin{equation}\label{farawaypoint2}
 \|u(x)\|, \|u'(x)\|\geq \frac{1}{2}\qquad\text{for any }x\in B_1(x_3).
 \end{equation}
 
 The proof of Theorem~\ref{mainthm} requires several steps. The first one is to estimate $\|u-u'\|$ on $\Gamma_{s}$,
 for $0<s\leq s_0$. This is obtained by a classical quantitative unique continuation.
 
 \begin{lem}\label{firstcontinuation}
Assume that $\varepsilon\leq 1/(2\rme)$. For any $0<s\leq s_0$, we have
\begin{equation}\label{stima1}
\|(u-u')(x)\|\leq \eta_{s}(\varepsilon)=E_1\varepsilon^{a^{m(s)}}\qquad\text{for any }x\in\Gamma_{s},
\end{equation}
with
\begin{equation}\label{numberest0}
m(s)\leq \frac{F_0}{s^N},
\end{equation}
where $E_1>0$, $F_0>0$ and $a$, $0<a<1$, are constants depending on the a priori data only.
 \end{lem}
  
\proof{.}
For any $x\in B_R\cap V_{s}$, $0<s\leq s_0$, let $\gamma$ be the curve connecting $x$ to $x_0$ as in the definition of $V_{s}$. Without loss of generality, we can assume that $\gamma$ is contained in $B_{R_1}$.

 We can construct a regular chain of balls, in the sense of \cite[Definition~5.1]{Ron08}, with respect to
 $B_{s}(\gamma)$ that from $x_0$ reaches $x$. The first ball is centred at $x_0$ and has radius less than or equal to 
 $\tilde{s}/8$. By Lemma~\ref{regullemma}, we have that
 $$\|u_p-u'_p\|_{L^{\infty}(B_{\tilde{s}/8}(x_0))},\ \|u_s-u'_s\|_{L^{\infty}(B_{\tilde{s}/8}(x_0))}\leq C\varepsilon.$$
 for a constant $C$ depending on the a priori data.
 
Then by a repeated use of the three-spheres inequality of Lemma~\ref{3sphereslemma} applied to $u_p-u'_p$ and $u_s-u'_s$ along this regular chain of balls, we obtain that
$$\|(u-u')(x)\|\leq \|(u_p-u'_p)(x)\|+\|(u_s-u'_s)(x)\| \leq E_1\varepsilon^{a^{m(s)}}
$$
where $E_1$ and $a$, $0<a<1$, depend on the a priori data only, and $m(s)$ denotes the number of times we have used the three-spheres inequality. It can be shown that $m(s)$ satisfies \eqref{numberest0}
for a constant $F_0$ depending on $R_1$ only.
Then we conclude the proof by using the continuity of $u$ and $u'$.\cvd

\smallskip
  
\begin{lem}\label{cruciallem}
Let $0\leq s\leq s_0$.
Assume that, for some $\eta$, $0<\eta\leq E$, we have
\begin{equation}
\|(u-u')(x)\|\leq \eta\qquad\text{for any }x\in\Gamma_{s}.
\end{equation}

Then
there exists a positive constant $\hat{C}_0$, depending on the a priori data only, such that
\begin{equation}\label{finalineq}
\|\nabla u\|^2_{L^2(A_{s})}\leq \hat{C}_0\hat{\eta}_s,
\end{equation}
where
\begin{equation}\label{tildeeta}
\hat{\eta}_s=(\eta+\tilde{C}_0s),
\end{equation}
$\tilde{C}_0$ as in \eqref{regularity}.
\end{lem}

\proof{.} We have
$\partial A_{s}\subset \partial K\cup \partial (D_{s})\cup\partial (D'_{s})$. We have that $u=0$ on $\partial K$ and, by \eqref{regularity}, $\|u\|\leq \tilde{C}_0s$ on $\partial (D_{s})$. By the same reasoning,
$\|u'\|\leq \tilde{C}_0s$ on $\partial (D'_{s})$. Hence
\begin{equation}\label{stima}
\|u(x)\|\leq \hat{\eta}_s=(\eta+\tilde{C}_0s)\qquad\text{for any }x\in\partial A_{s}.
\end{equation}
Since $u\in C^2(\overline{\Omega})$ and $A_{s}$ is a set of finite perimeter, an integration by parts leads to
$$
2\mu\int_{A_{s}}Eu : \overline{Eu}+\lambda\int_{A_{s}}\mathrm{div} (u) \overline{\mathrm{div} (u)} - \rho\omega^2\int_{A_{s}}u\cdot \overline{u}=\int_{\partial A_{s}}Tu\cdot\overline{u}.
$$
Then, by \eqref{regularity}, \eqref{Friedrichs}, \eqref{measurebdry+closeness} and \eqref{stima}, we have,
for a constant $M_1$ depending on the a priori data only,
\begin{multline}\label{firstineq}
M_1\hat{\eta}_s\geq
2\mu\int_{A_{s}}Eu : \overline{Eu}+\lambda\int_{A_{s}}\mathrm{div} (u) \overline{\mathrm{div} (u)} - \rho\omega^2\int_{A_{s}}u\cdot \overline{u}\\\geq \min\{2\mu,2\mu+\lambda\}\|Eu\|^2_{L^2(A_{s})} - 32\rho\omega^2 C(N)^2\tilde{H}_0^{\frac{1}{N}}
\left[\tilde{H}_0^{\frac{1}{N}}\|\nabla u\|^2_{L^2(A_{s})}+\|u\|^2_{L^2(\partial A_{s})}
\right]
\\
\geq \min\{2\mu,2\mu+\lambda\}\|Eu\|^2_{L^2(A_{s})} - 32\rho\omega^2 C(N)^2\left[\tilde{H}_0^{\frac{2}{N}}
\|\nabla u\|^2_{L^2(A_{s})}
+\tilde{H}_0^{\frac{1}{N}}(3C_4\hat{\eta}_s^2)\right].
\end{multline}
Then, by the idea of the proof of the first Korn inequality, we have
\begin{multline*}
\|Eu\|^2=\frac{1}{4}\sum_{i,j=1}^N\left|u^i_j+u^j_i\right|^2\\=\frac{1}{2}\sum_{i,j=1}^N\left|u^i_j\right|^2+
\frac{1}{4}\sum_{i,j=1}^N \left(u^i_j\overline{u}^j_i+\overline{u}^i_j u^j_i\right)=
\frac{1}{2}\sum_{i,j=1}^N\left|u^i_j\right|^2+
\frac{1}{2}\sum_{i,j=1}^N \left(u^i_j\overline{u}^j_i\right).
\end{multline*}
By the regularity of $u$ and of $A_{s}$, by two integrations by parts we obtain
\begin{multline*}
\int_{A_{s}}u^i_j\overline{u}^j_i=-\int_{A_{s}} u^i\overline{u}^j_{ij}+\int_{\partial A_{s}}u^i\overline{u}^j_i\nu^j=
-\int_{A_{s}} u^i\overline{u}^j_{ji}+\int_{\partial A_{s}}u^i\overline{u}^j_i\nu^j
\\=\int_{A_{s}} u^i_i\overline{u}^j_j-\int_{\partial A_{s}}u^i\overline{u}^j_j\nu^i
+\int_{\partial A_{s}}u^i\overline{u}^j_i\nu^j.
\end{multline*}
It follows that
$$\int_{A_{s}}\|Eu\|^2=\frac{1}{2}\int_{A_{s}}\|\nabla u\|^2+\frac{1}{2}\int_{A_{s}} \mathrm{div}(u)\overline{\mathrm{div}(u)}+
\frac{1}{2}\sum_{i,j=1}^N\left(-\int_{\partial A_{s}}u^i\overline{u}^j_j\nu^i
+\int_{\partial A_{s}}u^i\overline{u}^j_i\nu^j\right),
$$
thus
$$\int_{A_{s}}\|\nabla u\|^2\leq 2\int_{A_{s}}\|Eu\|^2+\sum_{i,j=1}^N\left(\int_{\partial A_{s}}u^i\overline{u}^j_j\nu^i
-\int_{\partial A_{s}}u^i\overline{u}^j_i\nu^j\right).$$
We conclude by Theorem~\ref{regularityteo} and \eqref{stima} that 
\begin{equation}\label{crucialineq}
\int_{A_{s}}\|\nabla u\|^2\leq 2\int_{A_{s}}\|Eu\|^2+6C_4N^2\tilde{C}_0\hat{\eta}_s.
\end{equation}
Coupling \eqref{firstineq} and \eqref{crucialineq}, we obtain
\begin{multline*}
M_1\hat{\eta}_s\geq \left(\frac{\min\{2\mu,2\mu+\lambda\}}{2}-32\rho\omega^2 C(N)^2\tilde{H}_0^{\frac{2}{N}}\right)\|\nabla u\|^2_{L^2(A_{s})}\\
-32\rho\omega^2 C(N)^2\tilde{H}_0^{\frac{1}{N}}(3C_4\hat{\eta}_s^2)-\min\{2\mu,2\mu+\lambda\}3C_4N^2\tilde{C}_0\hat{\eta}_s.
\end{multline*}
Since $\hat{\eta}_s$ is bounded by a constant depending on the a priori data only and, by \eqref{closeness} and\eqref{measurebdry+closeness}, we have that $\tilde{H}_0<H_1$, we can easily conclude the proof.\cvd

\smallskip

\begin{cor}\label{cor}
Under the same assumptions of Lemma~\textnormal{\ref{cruciallem}}, 
there exists a positive constant $\hat{C}_1$, depending on the a priori data only, such that
\begin{equation}\label{finalineq2}
\|u_p\|_{L^{\infty}(B_{c_3 d/4}(x_2))},\ \|u_s\|_{L^{\infty}(B_{c_3 d/4}(x_2))}\leq \frac{\hat{C}_1}{d^{(N+2)/2}}\hat{\eta}_s^{1/2},
\end{equation}
$\hat{\eta}_s$ as in \eqref{tildeeta}.
\end{cor}

\proof{.} Immediate by using \eqref{regestimate-b} and \eqref{finalineq}.\cvd

\smallskip

\begin{oss}\label{firststeposs}
Assume that $\varepsilon\leq \rme^{-\rme}/2$. Let $\hat{\eta}_s(\varepsilon)=(\eta_s(\varepsilon)+\tilde{C}_0s)$ as in \eqref{tildeeta} and with $\eta_s(\varepsilon)$ as in \eqref{stima1}. Then, by taking the minimum as $s$ varies in $(0,s_0]$, an easy computation shows that there exist positive constants
$\tilde{\varepsilon}$, $\tilde{\varepsilon}\leq \rme^{-\rme}/2$, and $\hat{C}_2$, depending on the a priori data only, such that if $0<\varepsilon \leq\tilde{\varepsilon}$ we have
\begin{equation}\label{finalineq3}
\|u_p\|_{L^{\infty}(B_{c_3 d/4}(x_2))},\ \|u_s\|_{L^{\infty}(B_{c_3 d/4}(x_2))}\leq \frac{\hat{C}_1}{d^{(N+2)/2}}\hat{\eta}(\varepsilon)^{1/2},
 \end{equation}
where
\begin{equation}\label{finalineq4}
\hat{\eta}(\varepsilon)= \hat{C}_2\left[\log(\log(1/\varepsilon))\right]^{-1/N}.
\end{equation}
In fact, let
$$s(\varepsilon)=\hat{C}_3\left[\log(\log(1/\varepsilon))\right]^{-1/N}$$
with $\hat{C}_3$ such that
$$\frac{\log(1/a)F_0}{\hat{C}_3^N}\leq 1/2.$$
Then
$$\varepsilon^{a^{m(s(\varepsilon))}}\leq \exp\left(-(\log(1/\varepsilon))^{1/2}\right).$$
It is enough to choose $\tilde{\varepsilon}$ such that for any $0<\varepsilon\leq\tilde{\varepsilon}$ we have $s(\varepsilon)\leq s_0$ and
$$\exp\left(-(\log(1/\varepsilon))^{1/2}\right)\leq \left[\log(\log(1/\varepsilon))\right]^{-1/N}.$$
\end{oss}

\smallskip

\begin{lem}\label{gettingoutlemma}
Let us assume that for some $d_1$, $0<d_1\leq c_3R/2$, there exists $x$ such that $B_{d_1}(x)\subset\Omega$ and,
for some $\hat{\eta}\leq \rme^{-\rme}/2$,
\begin{equation}\label{finalineq5}
\|u_p\|_{L^{\infty}(B_{d_1}(x))},\ \|u_s\|_{L^{\infty}(B_{d_1}(x))}\leq \hat{\eta},
\end{equation}

Then there exist a constant $\delta$, $0<\delta\leq \rme^{-\rme}/2$, and a positive constant $\hat{C}_4$, depending on the a
priori data only, such that if $0<\hat{\eta}\leq \delta$ we have
\begin{equation}\label{d1estimate}
d_1\leq \varphi(\hat{\eta})
\end{equation}
where
\begin{equation}\label{varphiestimate}
\varphi(\hat{\eta})\leq 2\rme R\left(\log(1/\hat{\eta})\right)^{-\hat{C}_4}.
\end{equation}
\end{lem}

\proof{.} We shall apply this lemma to $x=x_2$ and $d_1=c_3d/4$, see Remark~\ref{anotheross}.

By the uniform interior cone property of $\Omega$, we can find a direction $v$ such that
$\mathcal{C}=\mathcal{C}(x-(d_1)v,v,r_1,\theta_1)\subset \Omega$, for some $r_1>0$ and $\theta_1$, $0<\theta_1<\pi/2$, depending on $r$ and $L$ only. Moreover, we can find $\tilde{x}=x+sv$, for some $s\geq 0$, and $c_4$, $0<c_4<1$, depending on $r$ and $L$ only, such that $B_{c_4d_1}(x)$ and $B_{c_4r_1}(\tilde{x})$ are both contained in 
$\mathcal{C}$.

We can find $\gamma$, a piecewise $C^1$ curve, connecting $\tilde{x}$ with $x_3$ such that 
$B_{c_5r_1}(\gamma)\subset\Omega$, for some $c_5$, $0<c_5<1$, depending on $r$, $L$ and $R$ only. Without loss of generality, we can assume that $\gamma\subset B_{R_2}$ as well.

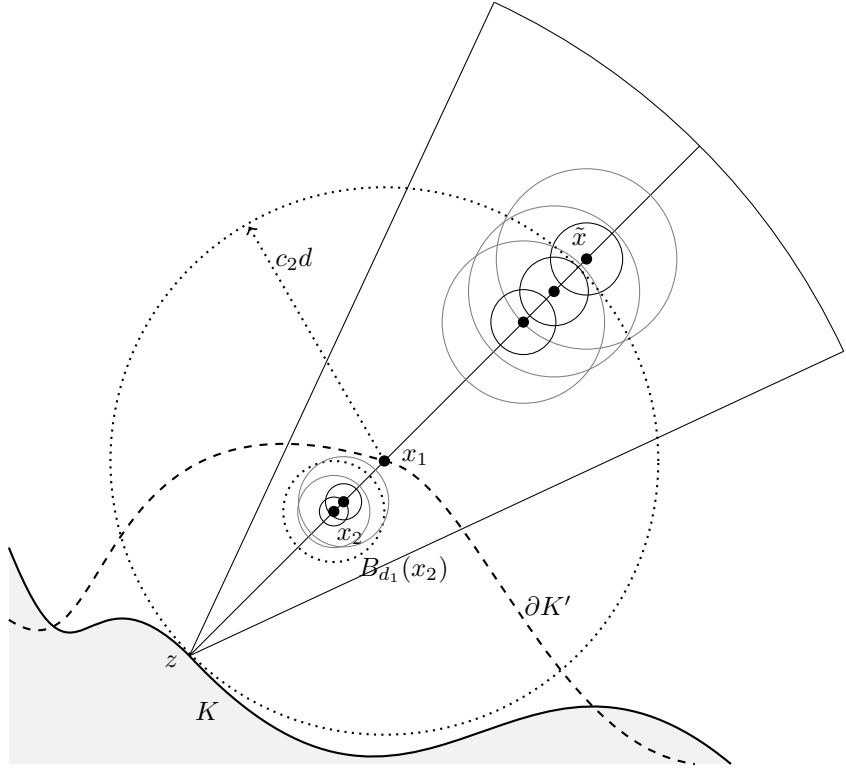
\begin{figure}[htb]
\centering
\begin{tikzpicture}[scale=4.8]
\draw[color=black] (0,0) to (.844,1.812);
\draw[color=black] (0,0) to (1.414,1.414);
\draw[color=black] (0,0) to (1.812,.844);
\draw [color=black] (1.812,.844) arc [radius=2, start angle=25, end angle=65];

\filldraw[color=black] (1.1,1.1) circle (.4pt);
\draw[color=gray](1.1,1.1) circle (0.25);
\draw[color=black](1.1,1.1) circle (0.1);

\filldraw[color=black] (1.01,1.01) circle (.4pt);
\draw[color=gray](1.01,1.01) circle (0.237);
\draw[color=black](1.01,1.01) circle (0.095);

\filldraw[color=black] (0.925,0.925) circle (.4pt);
\draw[color=gray](0.925,0.925) circle (0.225);
\draw[color=black](0.925,0.925) circle (0.09);


\draw[color=black, dotted, thick](.54,.54) circle (1.405*.54);
\filldraw[color=black] (0.54,0.54) circle (.4pt);
\draw[color=black, dotted, thick, ->] (.54,.54) to (.54-1.4105*.266,.54+1.4105*1.732*.266);
\node [right] at (.21, 1.1) {\small$c_2d$};
\draw[color=black, dotted, thick](.4,.4) circle (.14);

\node [right] at (.56, 0.55) {\small$x_1$};
\node [right] at (.38, 0.33) {\small$x_2$};

\node [right] at (.44, 0.24) {\small$B_{d_1}(x_2)$};

%
\filldraw[color=black] (0.427,0.427) circle (.4pt);
\draw[color=gray](0.427,0.427) circle (0.125);
\draw[color=black](0.427,0.427) circle (0.05);

\filldraw[color=black] (0.4,0.4) circle (.4pt);
\draw[color=gray](0.40,0.40) circle (0.1);
\draw[color=black](0.40,0.40) circle (0.04);

\path[fill=black!5] (-.5,0.3) .. controls (-.3,-.2) and (-.3,.3) .. (0,0) .. controls (.65,-.7) and (.9,.2) .. (1.5,-.3)
 .. controls (0,-.3) .. (-.5,-.3) .. controls (-.5,.0) .. (-.5,0.3);
\path[draw, color=black, thick] (-.5,0.3) .. controls (-.3,-.2) and (-.3,.3) .. (0,0) .. controls (.65,-.7) and (.9,.2) .. (1.5,-.3);
\node [below] at (0.05, -.1) {\small$K$};
\node [right] at (.9,.14) {\small$\partial K'$};
\node [below] at (-0.05, 0.03) {\small$z$};


\path[draw, color=black, dashed, thick] (-.5,.1) .. controls (-.2,-.1) and (-.3,.8) .. (0.54,0.54);
\path[draw, color=black, dashed, thick] (0.54,0.54) .. controls (.78,.5) and (1,-.25) .. (1.4,-.3);

\node [above] at (1.08, 1.11) {\small$\tilde{x}$};
\end{tikzpicture}
\caption{The geometric construction from $x=x_2$ to $\tilde{x}$.}\label{fig1}
\end{figure}

\begin{figure}[htb]
\centering
\begin{tikzpicture}[scale=2.8]
\draw[color=black] (0,0) to (.844,1.812);
\draw[color=black] (0,0) to (1.414,1.414);
\draw[color=black] (0,0) to (1.812,.844);
\draw [color=black] (1.812,.844) arc [radius=2, start angle=25, end angle=65];
\path[draw, color=black] (1.1,1.1) .. controls (2,.9) and (3,1.8)  .. (4,1.5);


\draw[color=black, dotted, thick](.54,.54) circle (1.405*.54);
\filldraw[color=black] (0.54,0.54) circle (.4pt);
\filldraw[color=black] (0.4,0.4) circle (.4pt);
\draw[color=black, dotted, thick, ->] (.54,.54) to (.54-1.4105*.266,.54+1.4105*1.732*.266);
\node [right] at (.21, 1.1) {\small$c_2d$};

\node [right] at (.56, 0.55) {\small$x_1$};
\node [right] at (.38, 0.33) {\small$x_2$};

\draw[color=black, dotted, thick](.4,.4) circle (.14);


\filldraw[color=black] (1.1,1.1) circle (.4pt);
\draw[color=gray](1.1,1.1) circle (0.25);
\draw[color=black](1.1,1.1) circle (0.1);

\filldraw[color=black] (1.23,1.08) circle (.4pt);
\draw[color=gray](1.23,1.08) circle (0.25);
\draw[color=black](1.23,1.08) circle (0.1);

\filldraw[color=black] (4,1.5) circle (.4pt);
\draw[color=gray](4,1.5) circle (0.25);
\draw[color=black](4,1.5) circle (0.1);

\filldraw[color=black] (3,1.487) circle (.3pt);
\draw[color=gray](3,1.487) circle (0.25);
\draw[color=black](3,1.487) circle (0.1);

\filldraw[color=black] (2.865,1.451) circle (.3pt);
\draw[color=gray](2.865,1.451) circle (0.25);
\draw[color=black](2.865,1.451) circle (0.1);

\filldraw[color=black] (2.732,1.413) circle (.3pt);
\draw[color=gray](2.732,1.413) circle (0.25);
\draw[color=black](2.732,1.413) circle (0.1);

\filldraw[color=black] (3.88,1.53) circle (.4pt);
\draw[color=gray](3.88,1.53) circle (0.25);
\draw[color=black](3.88,1.53) circle (0.1);

\path[fill=black!5] (-.5,0.3) .. controls (-.3,-.2) and (-.3,.3) .. (0,0) .. controls (.65,-.7) and (.9,.2) .. (1.5,-.3)
 .. controls (0,-.3) .. (-.5,-.3) .. controls (-.5,.0) .. (-.5,0.3);
\path[draw, color=black, thick] (-.5,0.3) .. controls (-.3,-.2) and (-.3,.3) .. (0,0) .. controls (.65,-.7) and (.9,.2) .. (1.5,-.3);
\node [below] at (0.05, -.1) {\small$K$};
\node [right] at (.9,.14) {\small$\partial K'$};
\node [below] at (-0.05, 0.03) {\small$z$};

\node [above] at (1.08, 1.11) {\small$\tilde{x}$};
\node [right] at (3.98, 1.47) {\small$x_3$};

\node [right] at (2, 1.27) {\small$\gamma$};


\path[draw, color=black, dashed, thick] (-.5,.1) .. controls (-.2,-.1) and (-.3,.8) .. (0.54,0.54);
\path[draw, color=black, dashed, thick] (0.54,0.54) .. controls (.78,.5) and (1,-.25) .. (1.4,-.3);

%
%
%
\end{tikzpicture}
\caption{From $\tilde{x}$ to $x_3$.}\label{fig2}
\end{figure}
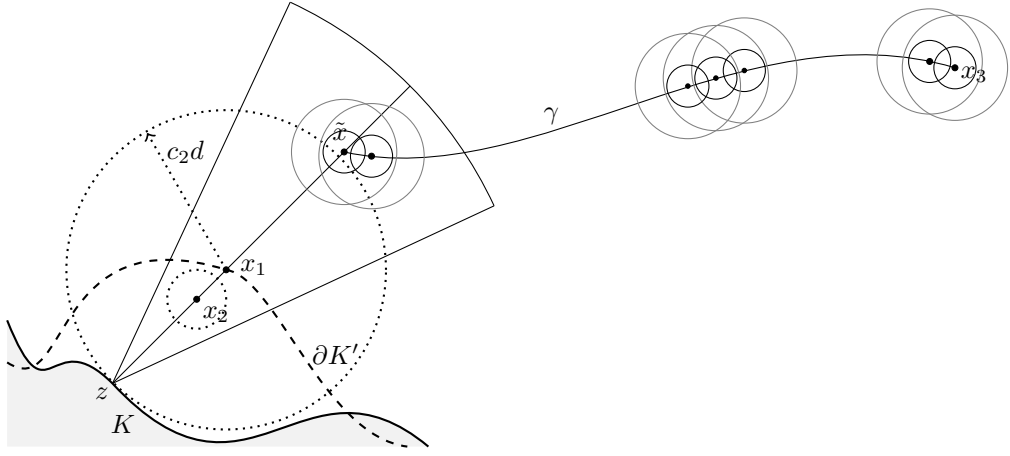

We can construct a regular chain of balls, again in the sense of \cite[Definition~5.1]{Ron08}, with respect to
$\mathcal{C}\cup B_{c_5r_1}(\gamma)$
that from $x$ reaches $x_3$, see the Step~I of the proof of \cite[Theorem~4.1]{Ron08} for details on this geometric construction which we just sketch now.
The first ball is centred at $x$ and has radius less than or equal to 
 $c_4d_1/8$. Then we proceed along the bisecting line of $\mathcal{C}$ until we reach 
 $\tilde{x}$. The construction from $x$ to $\tilde{x}$ is illustrated in Figure~\ref{fig1}. From $\tilde{x}$ to $x_3$ we proceed along the curve $\gamma$, see Figure~\ref{fig2} for an illustration.

Again by a repeated use of the three-spheres inequality of Lemma~\ref{3sphereslemma} applied to $u_p$ and $u_s$ along this regular chain of balls, we obtain that
$$\|u(x_3)\|\leq \|u_p(x_3)\|+\|u_s(x_3)\| \leq E_2\hat{\eta}^{a^n}
$$
where $E_2\geq 2\rme$ and $a$, $0<a<1$, depend on the a priori data only, and $n$ denotes the number of times we have used the three-spheres inequality. We can estimate $n$ as follows
\begin{equation}\label{numberest}
n\leq F_1\log(2\rme R/d_1)
\end{equation}
for some constant $F_1$ depending on the a priori data only.

It follows by \eqref{farawaypoint2} that
$$\frac12\leq E_2\hat{\eta}^{a^n},$$
consequently
$$-\log(2E_2)\leq -\log(1/\hat{\eta})a^n,$$
that is,
$$a^n\leq \frac{\log(2E_2)}{\log(1/\hat{\eta})}.$$
So, by \eqref{numberest},
$$\log(1/a)F_1\log(2\rme R/d_1)\geq n\log(1/a)\geq \log(\log(1/\hat{\eta}))-\log(\log(2E_2)).$$
The proof can now be concluded by an elementary computation.\cvd

\smallskip

\begin{oss}\label{anotheross}
Let us pick, in Lemma~\ref{gettingoutlemma}, $x=x_2$, $d_1=c_3 d/4$ and
$$\hat{\eta}=\frac{\hat{C}_1}{d^{(N+2)/2}}\hat{h}^{1/2}$$
with $\hat{h}=\hat{\eta}_s$ as in Corollary~\ref{cor} or $\hat{h}=\hat{\eta}(\varepsilon)$ as in Remark~\ref{firststeposs},
and assume that $\hat{h}\leq 1/(2\rme)$.
We conclude that either $\hat{\eta}\geq \delta$, that is,
\begin{equation}\label{caso1}
d^{(N+2)/2}\leq \frac{\hat{C}_1}{\delta} \hat{h}^{1/2},
\end{equation}
or $\hat{\eta}\leq \delta$ and, by \eqref{d1estimate} and \eqref{varphiestimate},
$$d\leq \frac{8\rme R}{c_3}\left(\frac{N+2}{2}\log(d)-
\log(\hat{C}_1)+\frac{1}{2}\log\left(\frac1{\hat{h}}\right)\right)^{-\hat{C}_4}.
$$
In this case we have that
$$\left(\frac{8\rme R}{c_3}\right)^{1/\hat{C}_4}\left(\frac{1}{d}\right)^{1/\hat{C}_4}\geq \frac{N+2}{2}\log(d)-
\log(\hat{C}_1)+\frac{1}{2}\log\left(\frac1{\hat{h}}\right),$$
that is,
$$\left(\frac{8\rme R}{c_3}\right)^{1/\hat{C}_4}\left(\frac{1}{d}\right)^{1/\hat{C}_4}+\frac{N+2}{2}\log\left(\frac1d\right)+
\log(\hat{C}_1)\geq \frac{1}{2}\log\left(\frac1{\hat{h}}\right).
$$
Since $d\leq 2R$, therefore $(1/d)\geq 1/(2R)$, we can find a positive constant $\hat{C}_5$, depending on the a priori data only, such that for any $0<d\leq 2R$ we have
$$\left(\frac{8\rme R}{c_3}\right)^{1/\hat{C}_4}\left(\frac{1}{d}\right)^{1/\hat{C}_4}+\frac{N+2}{2}\log\left(\frac1d\right)+
\log(\hat{C}_1)\leq \frac{1}{2}\left(\frac{\hat{C}_5}{d}\right)^{1/\hat{C}_4},$$
therefore
\begin{equation}\label{caso2}
d\leq \hat{C}_5\left(\log\left(\frac1{\hat{h}}\right)\right)^{-\hat{C}_4}.
\end{equation}
We can also find a positive constant $\hat{C}_6$, depending on the a priori data only, such that
for any $\hat{h}$, $0<\hat{h}\leq 1/(2\rme)$, we have
$$\frac{\hat{C}_1}{\delta} \hat{h}^{1/2}\leq \hat{C}_6\left(\log\left(\frac1{\hat{h}}\right)\right)^{-\hat{C}_4}.$$
Coupling \eqref{caso1} and \eqref{caso2}, we conclude that
\begin{equation}\label{final}
d\leq \hat{C}_7\left(\log\left(\frac1{\hat{h}}\right)\right)^{-\hat{C}_4},
\end{equation}
where $\hat{C}_7=\max\{\hat{C}_5,\hat{C}_6\}$.
\end{oss}

\smallskip

\proof{ of Theorem~\textnormal{\ref{mainthm}}.}
We can find $\hat{\varepsilon}$, depending on the a priori data only, such that $0<\hat{\varepsilon}\leq \tilde{\varepsilon}$ and for any $\varepsilon$,
$0<\varepsilon\leq \hat{\varepsilon}$, we have
$\hat{\eta}(\varepsilon)\leq 1/(2\rme)$, where $\hat{\eta}(\varepsilon)$ is defined in 
\eqref{finalineq4}. By the reasoning used in Remark~\ref{anotheross} with $\hat{h}=\hat{\eta}(\varepsilon)$, we deduce that
\begin{equation}\label{finalbis}
\tilde{d}\leq C_3d\leq C_3\hat{C}_7\left(\log\left(\frac1{\hat{\eta}(\varepsilon)}\right)\right)^{-\hat{C}_4},
\end{equation}
where we used \eqref{distancesoppositebis}. We have already obtained a quantitative estimate, which we can improve as follows.

Up to taking a smaller $\hat{\varepsilon}>0$, still depending on the a priori data only, for any $\varepsilon$,
$0<\varepsilon\leq \hat{\varepsilon}$, we have that
$\tilde{d}\leq \tilde{d}_0$ where $\tilde{d}_0$ is the constant of Lemma~\ref{RLGlemma}.
Then we can improve our estimate with a by now classical technique, which we sketch now.

Let us assume that $\varepsilon$,
$0<\varepsilon\leq \hat{\varepsilon}$, so that $\tilde{d}\leq \tilde{d}_0$. By Lemma~\ref{RLGlemma},
we have that
$\Omega^{ext}$ satisfies a uniform interior cone property, with constants $\tilde{r}_0$ and $\tilde{\theta}_0$ depending on $r$ and $L$ only. Let $z\in \Gamma$ and let $\mathcal{C}=\mathcal{C}(z,v,\tilde{r}_0,\tilde{\theta}_0)\subset \Omega^{ext}$
for a suitable direction $v$. For any $s$, $0<s\leq \tilde{s}_0$, with $\tilde{s}_0\leq2$ small enough, let $x(s)=z+sv$. By a completely analogous construction to the one used in Lemma~\ref{gettingoutlemma}, just by reversing the chain of balls,
we connect $x_0$ to $x(s)$ with a suitable regular chain of balls contained in $\Omega^{ext}$. The construction is again illustrated in Figures~\ref{fig2} and \ref{fig1}, by replacing $x_3$ with $x_0$ and $\tilde{x}$ with $x$, and assuming that $z\in\Gamma$ and that the cone and $\gamma$ are contained in $\Omega^{ext}$.

The repeated use of the three-spheres inequality applied to $u_p-u'_p$ and $u_s-u'_s$ along this chain, from $x_0$ to $x(s)$, allows us to estimate
$$
\|(u-u')(x(s))\|\leq  E_3\varepsilon^{a^{l(s)}}
$$
where $E_3\geq 2\rme$, $0<a<1$ and
$$
l(s)\leq F_2\log(2\rme/s).
$$
As usual, $\tilde{s}_0$, $E_3$, $a$, and $F_2$ can be chosen as depending on the a priori data only.
We conclude that, for any $z\in\Gamma$,
$$\|(u-u')(z)\|\leq  E_3\varepsilon^{a^{l(s)}}+\tilde{C}_0s\quad\text{for any }0<s\leq\tilde{s}_0.$$
By reasoning as in Remark~\ref{firststeposs}, we estimate the minimum as $s$ varies in $(0,\tilde{s}_0]$.
Let us further assume, without loss of generality, that for any $\varepsilon$,
$0<\varepsilon\leq \hat{\varepsilon}$, we have 
$$s(\varepsilon)=\frac{2\rme}{\log(1/\varepsilon)^{1/\hat{C}_8}}\leq \tilde{s}_0,$$
with $\hat{C}_8$ such that
$$\frac{\log(1/a)F_2}{\hat{C}_8}\leq 1/2.$$
Then
$$\varepsilon^{a^{l(s(\varepsilon))}}\leq \exp\left(-(\log(1/\varepsilon))^{1/2}\right).$$
It is enough to choose $\hat{\varepsilon}$ such that for any $0<\varepsilon\leq\hat{\varepsilon}$ we have $s(\varepsilon)\leq \tilde{s}_0$ and, calling $\hat{\beta}=1/\hat{C}_8$,
$$\exp\left(-(\log(1/\varepsilon))^{1/2}\right)\leq \log(1/\varepsilon)^{-\hat{\beta}}.$$

Therefore
we can find $\hat{\varepsilon}$, $0<\hat{\varepsilon}\leq \rme^{-\rme}/2$ and a positive constant $\hat{\beta}_0$, both depending on the a priori data only, such that for any $\varepsilon$, $0<\varepsilon\leq\hat{\varepsilon}$,
and
for any $z\in \Gamma$ we have
$$\|(u-u')(z)\|\leq \log(1/\varepsilon)^{-\hat{\beta}_0}=\eta(\varepsilon)  \leq 1/(2\rme).
$$
By the same argument used to prove \eqref{finalbis}, if we replace $\hat{\eta}(\varepsilon)$ with $\eta(\varepsilon)$,
we conclude that the following stability result holds.

There exist positive constants $\hat{\varepsilon}$, $0<\hat{\varepsilon}\leq \rme^{-\rme}/2$, $C$ and $\beta$, depending on the a priori data only, such that for any $\varepsilon$, $0<\varepsilon\leq \hat{\varepsilon}$,if
$$\|u-u'\|_{L^{\infty}(B_{\tilde{s}}(x_0),\mathbb{C}^N)}\leq \varepsilon$$
then
\begin{equation}\label{finalter}
\tilde{d}\leq 
C\left(\log\left(\log(1/\varepsilon)\right)\right)^{-\beta},
\end{equation}
where
$$C=C_3\hat{C}_7\hat{\beta}_0^{-\hat{C}_4}\qquad\text{and}\qquad \beta=\hat{C}_4.$$
 
If we consider the far-field error $\varepsilon_0$ instead of the error $\varepsilon$, we can find $\hat{\varepsilon}_0$, $0<\hat{\varepsilon}_0\leq \rme^{-\rme}/2$, such that for any $\varepsilon_0$, $0<\varepsilon_0\leq \hat{\varepsilon}_0$, we have that $\tilde{\eta}(\varepsilon_0)\leq \hat{\varepsilon}$, hence by replacing $\varepsilon$ with $\tilde{\eta}(\varepsilon_0)$ in \eqref{finalter}, we conclude that
$$
\tilde{d}\leq 
\hat{C}\left(\log\left(\log(1/\varepsilon_0)\right)\right)^{-\beta},
$$
where $\hat{C}=4^{\beta}C$. Therefore the proof of the main theorem is concluded.\cvd

\appendix

\section{Proof of Theorem~\ref{regularityteo}.}

In this appendix we sketch the proof of the regularity of solutions to \eqref{scatteringpbmrigid}.

Let us first observe that it is enough to prove that there exist positive constants $r_1$ and $\tilde{C}$, depending on $r$, $L$, $R$, $\alpha$ and the coefficients only, such that
\begin{equation}\label{globalboundedness}
\|u\|_{L^{2}(\Omega\cap B_{R+1})}\leq \tilde{C}
\end{equation}
and
\begin{equation}\label{localregularity}
\|u\|_{C^{2}\left(\overline{\Omega}\cap\overline{B_{r_1}(z)}\right)}\leq \tilde{C}\qquad\text{for any }z\in\partial\Omega.
\end{equation}

Assuming we have proved \eqref{globalboundedness} and \eqref{localregularity}, we conclude
the proof of Theorem~\ref{regularityteo}. By the techniques developed in the proof of Lemma~\ref{regullemma}, we first show that by \eqref{globalboundedness} we have
\begin{equation}\label{localregularityter}
\|u\|_{C^{2}\left(\overline{B_{R+3/4}}\backslash B_{R+1/4}\right)}\leq \tilde{C}_2
\end{equation}
and then, using also \eqref{localregularity}, that
\begin{equation}\label{localregularitybis}
\|u\|_{C^{2}\left(\overline{\Omega}\cap\overline{B_{R+3/4}}\right)}\leq \tilde{C}_2
\end{equation}
for a constant $\tilde{C}_2$ depending on $r$, $L$, $R$, $\alpha$ and the coefficients only.

The estimate \eqref{localregularitybis} implies that $u_p$ and $u_s$, and thus $u_p^{scat}$ and $u_s^{scat}$,
are uniformly bounded in $\overline{\Omega}\cap\overline{B_{R+3/4}}$ by a constant depending on $\tilde{C}_2$ and the coefficients of \eqref{Navier} only.
By standard regularity estimates for solutions to the Helmholtz equation, not very different from what we used in the proof of
Lemma~\ref{regullemma}, we infer that
$\|u_p^{scat}\|+\|\nabla u_p^{scat} \nu\|$ is bounded on $\partial B_{R+1/2}$ by a constant depending on $\tilde{C}_2$, $R$, and the coefficients of \eqref{Navier} only, $\nu$ being the exterior normal to $B_{R+1/2}$. We point out that for this estimate \eqref{localregularityter}, thus \eqref{globalboundedness}, is enough. Then, for any $x$ with $\|x\|>R+1/2$ we have
$$u_p^{scat}(x)=\int_{\partial B_{R+1/2}}\left(\frac{\partial\phi_{\omega_p}(x,y)}{\partial\nu(y)}u_p^{scat}(y)-\phi_{\omega_p}(x,y)\nabla u_p^{scat}(y) \nu(y)\right)\, d\sigma(y).$$
Then, by the regularity and decay properties of $\phi_{\omega_p}$, it is not difficult to prove \eqref{decay}
for what concerns $u_p^{scat}$ and to bound
$\|u_p\|_{C^{2}\left(\Omega\backslash B_{R+3/4}\right)}$. A completely analogous argument applied to $u_s$ completes the proof.

The proof of \eqref{globalboundedness} and \eqref{localregularity} will be done in two steps. In the first step we prove
\eqref{globalboundedness}, in the second we prove \eqref{localregularity}.

\smallskip

\noindent
{\bf{Step I:}} The estimate \eqref{globalboundedness} is proved by a continuity argument which is inspired by Mosco convergence. We sketch the proof, for details we refer to \cite{Men-Ron} where the argument is fully developed in the acoustic case for the much harder Neumann boundary condition and for much more general classes of scatterers.

Let $\mathcal{A}$ be the class of obstacles contained in $\overline{B_R}$ whose interior is a Lipschitz open set with constants $r$ and $L$. In \cite[Section~2]{LPRX} it is proved that $\mathcal{A}$ is compact with respect to the Hausdorff distance. We claim that there exists a constant $\tilde{C}$, depending on $r$, $L$, $R$, and the coefficients only, such that \eqref{globalboundedness} holds for any $u$ solution to \eqref{scatteringpbmrigid} with $K\in \mathcal{A}$.

We argue by contradiction. Let us assume that there exists a sequence $\{K_n\}_{n\in\mathbb{N}}\subset\mathcal{A}$
such that, calling $u_n$ the solution to \eqref{scatteringpbmrigid} with $K$ replaced by $K_n$, we have for any $n\in\mathbb{N}$
$$\|u_n\|_{L^2(\Omega_n\cap B_{R+1})}=a_n\geq n.$$
We always extend $u_n$ to $0$ in $D_n$, the interior of $K_n$, so that $u_n\in H^1(B_{R+1})$.
Let $v_n=u_n/a_n$. We have that, with the usual extension to $0$ outside $\Omega_n$,
$$\|v_n\|_{L^2(B_{R+1})}=1\quad\text{for any }n\in\mathbb{N}.$$
By the same argument we used before, it is not difficult to show that
$\|v_n\|_{C^2(\overline{B_{R+3/4}}\backslash B_{R+1/4})}$ is bounded by a constant not depending on $n$. Let $\chi\in C^{\infty}_0(B_{R+1})$ be a cutoff function such that $0\leq \chi \leq 1$ in $B_{R+1}$ and $\chi=1$ in
$B_{R+1/2}$. We have that $w_n=v_n\chi\in H^1_0(B_{R+1}\backslash K_n)\subset H^1_0(B_{R+1})$ and $w_n$ solves
$$\left\{\begin{array}{ll}
\mu\Delta w_n+(\lambda+\mu)\nabla^T(\mathrm{div} (w_n))+\rho \omega^2 w_n=f_n &\text{in }B_{R+1}\backslash K_n\\
w_n=0&\text{on }\partial (B_{R+1}\backslash K_n)
\end{array}\right.
$$
with $\|f_n\|_{L^{\infty}(B_{R+1}\backslash K_n)}$ bounded by a constant not depending on $n$. By the weak formulation and first Korn inequality, we deduce that $\{w_n\}_{n\in\mathbb{N}}$ is bounded in $H^1(B_{R+1})$, hence
$\{v_n\}_{n\in\mathbb{N}}$ is bounded in $H^1(B_{R+1})$ as well.

Passing to subsequences, without loss of generality, we can assume that, as $n\to\infty$,
$K_n$ converges to $K\in\mathcal{A}$ in the Hausdorff distance, $v_n$ converges to $\tilde{v}$ weakly in $H^1(B_{R+1})$ and strongly in $L^2(B_{R+1})$. In particular,
\begin{equation}\label{contrad}
\|\tilde{v}\|_{L^2(B_{R+1})}=1.
\end{equation}

 One can easily show that $\tilde{v}= 0$ in a weak sense on $\partial K$ and that
$\tilde{v}$ solves \eqref{Navier} in $B_{R+1}\backslash K$. Since $v_n^{inc}=u^{inc}/a_n$, for any $n\in\mathbb{N}$,
we actually have that $v_n^{scat}$ converges to $\tilde{v}$ weakly in $H^1(B_{R+1})$ and strongly in $L^2(B_{R+1})$.
We infer that, possibly passing to a further subsequence by using a diagonal argument, $v_n^{scat}$ converges to a function $\tilde{v}$ in $L^2(B_r)$ for any $r>R$. Such a function $\tilde{v}$ is a radiating solution to Navier equation in $\Omega=\mathbb{R}^N\backslash K$ and $\tilde{v}=0$ on $\partial K$ in a weak sense. By uniqueness of the scattering problem we deduce that $\tilde{v}=0$ and this contradicts \eqref{contrad}.

\smallskip

\noindent
{\bf{Step II:}}
Let $z \in \partial \Omega$ and let $R_z$ be the unitary matrix transforming $\nu(z)$ to $(0,\ldots,1)$ where $\nu(z)$ is the exterior normal to $\partial \Omega$ at $z$. We call $V(x)=R_z u(R^{-1}_z(x))$ and we note that $V$ satisfies \eqref{Navier} in $R_z(\Omega)$. In other words, we can assume, without loss of generality, that there exists 
a
function $\phi_z:\mathbb{R}^{N-1}\to\mathbb{R}$, with $\|\phi_z\|_{C^{2,\alpha}(\mathbb{R}^{N-1})}\leq L$,
such that for any $y\in B_r(z)$ we have, without any further rigid transformation,
$$x=(x',x_N)\in\Omega\quad \text{if and only if}\quad  x_N<\phi_z(x').$$
We define
$\tilde{u}(\xi)=u(x)$ where $\xi=F(x)$ is defined as
\begin{equation}\label{Parametrisation}
\left \{
\begin{array}{l}
\xi'=x'\\
\xi_N=x_N-\phi_z(x').
\end{array} \right.
\end{equation}

Here, and in the sequel, for any $s>0$ we denote $\Sigma_{s}(z)=B_{s}(z)\cap\{\xi_N> z_N\}$ and
 $\Gamma_s(z)=B_{s}(z)\cap\{\xi_N= z_N\}$.
By the regularity properties of $\phi_z$, we can infer that, for a positive constant $r_2$, depending on $r$ and $L$ only, we have that $\tilde{u}$ satisfies, for any $i=1,\ldots,N$,
\begin{multline*}
\mu\left[\Delta_{\xi}\tilde{u}^i(\xi)-2\nabla_{\xi}\tilde{u}^i_N(\xi)\cdot \nabla_x\phi(x)+\tilde{u}^i_{NN}(\xi)\|\nabla_x \phi(x)\|^2\right]\\
+(\lambda+\mu)\left[(\mathrm{div}_{\xi}\tilde{u})_i(\xi)-(\mathrm{div}_{\xi}\tilde{u})_N(\xi)\phi_i(x)+\tilde{u}_{Ni}(\xi)\cdot\nabla_x \phi(x)-(\tilde{u}_{NN}(\xi)\cdot\nabla_x\phi(x))\phi_i(x)
\right]\\-\mu\tilde{u}^i_N(\xi)\Delta_x\phi(x)+(\lambda+\mu)\tilde{u}_N(\xi)\cdot\nabla_x(\phi_i)(x)+\rho\omega^2\tilde{u}^i(\xi)=0\quad\text{in }\Sigma_{r_2}(z)
\end{multline*}
and
$$\tilde{u}=0\quad\text{on }\Gamma_{r_2}(z).$$
Here $\phi:\mathbb{R}^N\to\mathbb{R}$ is defined as $\phi(x)=\phi_z(x')$ for any $x\in\mathbb{R}^N$ and, in the formula,
$x=F^{-1}(\xi)$ everywhere. By Remark~\ref{normaloss}, we have that $\nabla_x\phi(z)=0$. Hence, in a suitable neighbourhood of $z$, the principal part of the second order system solved by $\tilde{u}$ is a small perturbation of the Lam\'e system
$\mu\Delta\tilde{u}+(\lambda+\mu)\nabla^T(\mathrm{div}(\tilde{u}))$. We can conclude that it is elliptic and satisfies all the conditions of \cite[Chapter~1]{ADN}, as well as the boundary condition $\tilde{u}=0$ on $\Gamma_{r_2}(z)$
satisfies all the conditions of \cite[Chapter~2]{ADN}. Moreover, the coefficients of the elliptic system are bounded in $C^{0,\alpha}$.

As an intermediate step, we show that there exist $r_3>0$ and $\tilde{C}_3$, depending on $r$, $L$, $R$, $\alpha$ and the coefficients only, such that 
\begin{equation}\label{localboundedness}
\|\tilde{u}\|_{L^{\infty}\left(\Sigma_{r_3}(z)\right)}\leq \tilde{C}_3.
\end{equation}

By using \cite[Theorem~10.4]{ADN} with a suitable cutoff function,
for any $s$, $0<s\leq r_2$, and any real $p$, $p\geq 2$, we can find $s_1$, $0<s_1<s$, and $C_1$, depending on $s$, $p$ and
$r$, $L$, $R$, $\alpha$ and the coefficients only, such that
\begin{equation}\label{W2p}
\|\tilde{u}\|_{W^{2,p}(\Sigma_{s_1}(z))}\leq C_1\|\tilde{u}\|_{L^p(\Sigma_{s}(z))}.
\end{equation}
By \eqref{globalboundedness}, we can control $\|\tilde{u}\|_{W^{2,2}(\Sigma_{\tilde{s}}(z))}$ for a suitable positive constant $\tilde{s}$. By Sobolev inequality, we infer that  $\tilde{u}$ belongs to $L^{p_1}(\Sigma_{\tilde{s}}(z))$ for some $p_1>2$, thus, repeating the argument, $\tilde{u}$ belongs to $W^{2,p_1}(\Sigma_{\tilde{s}_1}(z))$ for a smaller positive constant $\tilde{s}_1$. With a bootstrap argument, after a finite number $m$ of steps, which depends on $N$ only, we obtain that $\tilde{u}$ belongs to $W^{2,p_m}(\Sigma_{\tilde{s}_m}(z))$ for a positive constant $\tilde{s}_m$ and $p_m>N$. By a final application of Sobolev inequality, we conclude that \eqref{localboundedness} holds.

Once \eqref{localboundedness} is established, by the standard estimates of \cite[Theorem~9.2]{ADN}, we can control the
$C^{2,\alpha}$-norm of $\tilde{u}$ in $\overline{\Sigma_{r_3/16}(z)}$. Going back to the usual coordinates, \eqref{localregularity} can be finally proved.

\end{document}